\documentclass{amsart}
\usepackage{epsfig}
\usepackage{amsthm}
\usepackage{amssymb}
\usepackage{amsmath}
\usepackage{amscd}

\newtheorem {Theorem}  {Theorem}
\numberwithin{Theorem}{section}
\numberwithin{equation}{section}
\newtheorem {Lemma}[Theorem]  {Lemma}

\theoremstyle{definition}
\newtheorem{Definition}[Theorem]{Definition}
\theoremstyle{remark}
\newtheorem{Remark}[Theorem]{Remark}

\newtheorem {Corollary}[Theorem]{Corollary}

\expandafter\chardef\csname pre amssym.def at\endcsname=\the\catcode`\@ \catcode`\@=11
\def\undefine#1{\let#1\undefined}
\def\newsymbol#1#2#3#4#5{\let\next@\relax
 \ifnum#2=\@ne\let\next@\msafam@\else
 \ifnum#2=\tw@\let\next@\msbfam@\fi\fi
 \mathchardef#1="#3\next@#4#5}
\def\mathhexbox@#1#2#3{\relax
 \ifmmode\mathpalette{}{\m@th\mathchar"#1#2#3}%
 \else\leavevmode\hbox{$\m@th\mathchar"#1#2#3$}\fi}
\def\hexnumber@#1{\ifcase#1 0\or 1\or 2\or 3\or 4\or 5\or 6\or 7\or 8\or
 9\or A\or B\or C\or D\or E\or F\fi}
\font\teneufm=eufm10 \font\seveneufm=eufm7 \font\fiveeufm=eufm5
\newfam\eufmfam
\textfont\eufmfam=\teneufm \scriptfont\eufmfam=\seveneufm
\scriptscriptfont\eufmfam=\fiveeufm

\catcode`\@=\csname pre amssym.def at\endcsname
\newcounter{remark}
\newcommand{\bg}{\begin{equation}}
\newcommand{\ed}{\end{equation}}
\newcommand{\bga}{\begin{eqnarray}}
\newcommand{\eda}{\end{eqnarray}}

\def\cbdu{\hfill{$\Box$}}

\renewcommand{\div}{\mbox{div}}

\def  \12  {{\frac{1}{2}}}

\begin{document}

\title[Decay and Existence for the Viscous Camassa-Holm Equations]{On Questions of Decay and Existence for the Viscous Camassa-Holm Equations}

\author[Clayton Bjorland ]{Clayton Bjorland \\Department of Mathematics\\ UC Santa Cruz\\ Santa Cruz, CA 95064, USA\\}


\author[Maria E. Schonbek ]{\\Maria E. Schonbek \\Department of Mathematics\\ UC Santa Cruz\\ Santa Cruz, CA 95064, USA}


\thanks{The work of M. Schonbek was partially supported by NSF Grant DMS-0600692.\\  The work of C. Bjorland was partially supported by NSF Grant OISE-0630623.}

\email{cbjorland@math.ucsc.edu} \email{schonbek@math.ucsc.edu}

\keywords{Navier-Stokes-alpha, Camassa-Holm, existence, regularity, decay}

\subjclass[2000]{35B40; 35K55} \

\date{\today}



\bigskip

\begin{abstract}
We consider the viscous $n$-dimensional Camassa-Holm equations, with $n=2,3,4$ in the whole space.  We establish
existence and regularity of the solutions and study the large time behavior of the solutions in several
Sobolev spaces.  We first show that if the data is only in $L^2$ then the solution decays without a rate and
that this is the best that can be expected for data in $L^2$.  For solutions with data in $H^m\cap L^1$ we obtain
decay at an algebraic rate which is optimal in the sense that it coincides with the rate of the underlying linear part.

\textbf{Quelques questions de decroissance et existence
pour les equations visqueuses de Camassa-Holm.}\\
R\'{e}sum\'{e} : On consid\`{e}re les \'{e}quations visqueuses de
Camassa--Holm dans $\mathbb{R}^n$, $n=2,3,4$. Nous
\'{e}tablissons l'existence
et regularit\'{e} des solutuions. Nous \'{e}tudions le
comportament asymptotique des solutions dans
plusieurs espaces de Sobolev quand le temps tend
vers l'infini. On montre que si la donn\'{e}e
est seulement dans $L^2$ la solution decro\^{i}t vers
zero, mais la decroissance ne peux \^{e}tre uniforme.
Pour les solutions avec de donn\'{e}e dans $L^1 \cap H^m$
on obtient  une decroissance alg\'{e}brique avec une
vitesse qui est optimale dans le sens que c'est la
m\^{e}me que pour les solutions correspondant a
l'\'{e}quation lin\'{e}aire.
\end{abstract}

\maketitle

\section{Introduction}
The Viscous Camassa-Holm equations (VCHE) are commonly written
\begin{align}\label{VCHE:PDE}
v_t +u\cdot\nabla v+v\cdot \nabla u^T +\nabla \pi &= \nu\triangle v \\
u-\alpha^2\triangle u&=v \notag \\
\nabla\cdot v&=0 \notag
\end{align}
Here we adopt the notation $(v\cdot \nabla u^T)_i = \sum_j v_j\partial_iu_j$.
These equations rose from work on shallow water
equations \cite{MR1234453}, which led to
\cite{MR1627802}, \cite{MR1853633}, where the equations are
derived by considering variational principles and
Lagrangian averaging.  In light of this derivation the equations are
sometimes called the Lagrangian Averaged Navier-Stokes equations. In
\cite{MR1837927}, the equations were derived as a ``filtered''
Navier-Stokes equation, which obeys a modified Kelvin circulation
theorem along filtered velocities.  In this setting they are
sometimes referred to as the Navier-Stokes-$\alpha$ equations, where
$\alpha$ is the parameter in the filter. Solutions to the VCHE are
closely related to solutions of the famous Navier-Stokes equation
(NSE), but the filter allows bounds that are currently unobtainable
for the NSE, making them in some ways better suited for computational turbulence study,
see \cite{htlans}.

In \cite{MR1837927}, \cite{MR1878243}
these equations were studied in relation to turbulence theory, this treatment includes existence and
uniqueness theorems on the torus in three dimensions.  The two
dimensional case was considered on the torus and the sphere in
\cite{MR2031580}.  Global existence
and uniqueness in three dimensions was proved on bounded domains
with zero (non-slip) boundary conditions in \cite{MR1853633}.  These equations have also been
studied in terms of large eddy simulation and turbulent pipe flow in
\cite{MR1745983},\cite{MR1721139},\cite{MR1719962},
and \cite{domholm2001}.  In this paper we extend the current existence theorems and study
the large time behavior of solutions.

This paper is organized as follows.  Section two consists of notation and conventions used throughout.  Section three contains preliminary discussion of the VCHE and several useful lemmas. In section four we state existence and uniqueness results for the VCHE, proofs of these statements are contained in the appendix.In the next two sections we continue the decay
program of M. E. Schonbek, \cite{MR1356749}, \cite{MR571048},
\cite{MR775190}, \cite{MR1312701},\cite{MR1396285}. The main result of chapter five considers solutions of the VCHE in the whole
space and we prove that the energy of a solution corresponding to data
only in $L^2(\mathbb{R}^n)$ decays to zero following the arguments
in \cite{MR1432588}.  We then demonstrate, by constructing counter examples, that no uniform rate of
decay can exists which depends only on the initial energy. In chapter six we consider decay for solutions with initial
initial data in $L^1\cap L^2$.  We show, using the
Fourier Splitting Method, that the energy of a solution decays at
the rate expected from the linear part, this is the same rate of
decay as solutions to the NSE.  For solutions with initial data in
$H^m\cap L^1$ we calculate the decay of derivatives using again the
Fourier Splitting Method with an inductive argument.  In section seven we examine how solutions of the VCHE approach solutions of the NSE strongly on intervals of regularity for the NSE.

\section{Notation}
In this paper, $L^p$ denotes the standard Lebesgue space with norm
$\|\phi\|_p=(\int |\phi|^p)^{1/p}$.  We use $<u,v>=\int uv$ to
denote the standard inner product on the Hilbert space $L^2$.
Compactly supported solenoidal vector fields (subsets of $\Sigma=\{\phi\in
C_0^\infty(\Omega)|\nabla\phi=0\}$) will be needed to describe
incompressible solutions with zero boundary conditions. $L_\sigma^p$ will denote the completion of
$\Sigma$ in the norm $\|\cdot\|_p$.  $W^{m,p}$ will be used to
denote the standard Sobolev spaces with the convention that
$H^m=W^{m,2}$ (and $L^2=H^0$).  The completion of $\Sigma$ under the $H^m$ norm will
be denoted by $H_\sigma^m$ and $(H^m_\sigma)'$ will be the dual
space. To denote the Fourier Transform of a function $\phi$ we will
use either $\hat{\phi}$ or $\mathcal{F}(\phi)$, with $\check{\phi}$
or $\mathcal{F}^{-1}(\phi)$ the inverse transform.  Throughout we will use $C$ to denote an arbitrary constant which may change line to line, to emphasis the dependence of a constant on a number, say $\nu$, we will write $C(\nu)$.

\section{Preliminaries}

The Kelvin-filtered Navier-Stokes equations (KFNSE) are given by the formula
\begin{align}
\frac{\partial v}{\partial t}+u\cdot\nabla v+ v\cdot \nabla u^T+\nabla \pi &=\nu\triangle v\notag \\
\nabla \cdot v=\nabla \cdot u&=0 \notag\\
v&=\mathcal{O}u\notag
\end{align}
In the above, $u=g\ast v$ represents a spatially filtered
fluid velocity and $\mathcal{O}$ is the inverse of this
convolution.  The term $u\cdot\nabla v$ is similar to ``mollifying''
the Navier-Stokes equations, originally done by Leray,
\cite{Leray}, to approximate solutions.  The term $v\cdot\nabla
u^T=\sum v_j\nabla u_j$ allows the solution to obey a modification
of the Kelvin circulation theorem where circulation is conserved
around a loop moving with the filtered velocity $u$.  In two and three dimensions, using the
identity
\begin{equation}\label{identity:vector}
u\cdot\nabla v +\sum v_j\nabla u_j= - u\times(\nabla\times v)+\nabla(v\cdot u)
\end{equation}
and including the term $\nabla(v\cdot u)$ in the pressure, the KFNSE
can be written as
\begin{align}
\frac{\partial v}{\partial t}+\nabla \pi &= u \times (\nabla \times v) + \nu\triangle v\notag\\
\nabla \cdot u=\nabla \cdot v&=0\notag\\
v&=\mathcal{O}u\notag
\end{align}
The following lemma will show that the bilinear term in the Kelvin-filtered
Navier-Stokes equations behaves similar to the bilinear term in the
Navier-Stokes equations.

\begin{Lemma}\label{KFNSE:bilinear}
Let $u$ and $v$ be smooth divergence free functions with compact support, then
\begin{align}
<u\cdot\nabla v,u>+<v\cdot \nabla u^T,u>=0 \notag\\
<u \times (\nabla \times v),u>=0 \notag
\end{align}
\end{Lemma}
\begin{proof}
The second equality is a consequence of the first, the identity
(\ref{identity:vector}), and the fact that $u$ is divergence free.
To see the first inequality we just need to rearrange the terms and
then integrate by parts
\begin{equation}
\sum_{i,j}\int_{\mathbb{R}^n} v_j\partial_i u_j u_i\,dx =
-\sum_{i,j}\int_{\mathbb{R}^n} u_i\partial_i v_j u_j\,dx\notag
\end{equation}
\end{proof}
Using this lemma, we can formally multiply the KFNSE by $u$ to find
\begin{equation}\label{KFNSENERGY}
<\frac{\partial}{\partial t} v,u>+ \nu<\nabla v,\nabla u>=0
\end{equation}

By choosing $\mathcal{O}$ to be the Helmholtz operator
$\mathcal{O}=1-\alpha^2\triangle$ we recover the Viscous
Camassa-Holm equations
\begin{align}
v_t +u\cdot\nabla v+v\cdot \nabla u^T +\nabla \pi &= \nu\triangle v\notag\\
u-\alpha^2\triangle u&=v\notag\\
\nabla\cdot v&=0\notag
\end{align}
In the case of the VCHE, (\ref{KFNSENERGY}) becomes
\begin{equation}\label{VCHE:energy}
\frac{1}{2}\frac{d}{d t}(<u,u>+\alpha^2<\nabla u, \nabla u>)+\nu(<\nabla u,\nabla u>+\alpha^2<\triangle u,\triangle u>)=0
\end{equation}
This relation gives a priori estimates on $u$:
\begin{align}\label{VCHE:aprioriestimate}
\|u(\cdot,t)\|^2_2+\alpha^2\|\nabla u(\cdot,t)\|^2_2 +2\nu \int^t_0\|\nabla u(\cdot,t)\|_2^2\,dt +2\nu \alpha^2 \int^t_0\|\nabla^2 u(\cdot,t)\|_2^2\,dt \notag\\
\leq \|u_0\|^2_2 +\alpha^2\|\nabla u_0\|_2^2
\end{align}

\section{Existence of Solutions for the VCHE}

Existence and uniqueness of solutions for the VCHE
on periodic domains in three dimensions was proved first in \cite{MR1878243} using
the Galerkin method.  The most general existence and uniqueness
theorems in three dimensions are provided in \cite{MR1853633} which
relies on a fixed point argument.   The theorems in \cite{MR1853633}
assume the initial data $u_0\in H^1_0\cap H^s$ with $s\in[3,5)$ and
$u=Au=0$ on the boundary, where $A$ is the Stokes operator. Here we state
extended results which cover the whole
space in dimensions $2\leq n\leq 4$, proofs are included in the appendix.  As an intermediate step, we
provide a new existence proof on bounded domains in dimensions
$2\leq n\leq 4$, using the Galerkin Method, with initial data
$v_0\in L^2$, and $u=v=0$ on the boundary.  Our bounded result in
three dimensions is slightly stronger then \cite{MR1853633}, by
assuming $v_0\in L^2$ we have only implied $u\in H_\sigma^2$.

\begin{Definition}\label{VCHE:weakdefn1}
Let $\Omega\subset \mathbb{R}^n$ be any open bounded subset or $\Omega=\mathbb{R}^n$, $n=2,3,4$.  A weak solution to the VCHE (\ref{VCHE:PDE}), with zero (no-slip)
boundary conditions in the case of $\Omega$ bounded, is a pair of functions, $u$, $v$, such that
\begin{align}
v&\in L^\infty ([0,T];L_\sigma^2(\Omega))\cap L^2 ([0,T];H_\sigma^1(\Omega))\notag\\
\partial_t v &\in L^2([0,T];(H_\sigma^1)'(\Omega))\notag\\
u&\in L^\infty ([0,T];H^2_\sigma(\Omega))\cap L^2 ([0,T];H_\sigma^3(\Omega))\notag
\end{align}
as well as $v(x,0)=v_0$, and for any $\phi \in L^2([0,T];H_\sigma^1(\Omega))$
with $\phi(T)=0$ the following equalities are satisfied:
\begin{align}
-\int_0^T<v,\partial_t\phi>\,ds &+ \int_0^T<u\cdot\nabla v,\phi>\,ds \notag\\
 &+ \int_0^T<\phi\cdot\nabla u,v>\,ds +\nu\int_0^T<\nabla v,\nabla \phi>\,ds =<v_0,\phi(0)> \notag
\end{align}
and for a.e. $t\in[0,T]$,
\begin{align}
<u,\phi> +\alpha^2<\nabla u,\nabla \phi> = <v,\phi>\notag
\end{align}
\end{Definition}

\begin{Theorem}\label{VCHE:existence}
Let $\Omega\subset\mathbb{R}^n$ be an open bounded set with smooth boundary
or $\Omega=\mathbb{R}^n$, $n=2,3,4$.
Given initial data $v_0\in H^M_\sigma(\Omega)$, $M\geq 0$, 
there exists a unique weak solution to the VCHE (\ref{VCHE:PDE}) in the sense of Definition \ref{VCHE:weakdefn1}.
This solutions satisfies the estimate (\ref{VCHE:aprioriestimate}) as well as
\begin{equation}\label{VCHE:regularityestimate}
\|\partial^p_t \nabla^{m}v(t)\|_2^2 +\nu\int_0^t\|\partial_t^p \nabla^{m+1}v(s)\|_2^2\,ds \leq C(\|v_0\|_{H^M_0})
\end{equation}
for all $m+2p\leq M$.
\end{Theorem}
\begin{proof}
Existence is given by Theorems \ref{VCHE:boundedexistence} and
\ref{VCHE:unboundedexistence} in the appendix.  The regularity
statement is Theorem \ref{VCHE:regularity} and uniqueness is Theorem
\ref{VCHE:uniequnesstheorem}.  The proofs follow from the construction of approximate solutions
using the Galerkin method on bounded domains.  A priori bounds are obtained through energy methods.  Using a compactness lemma we are
able to find a strongly convergent subsequence which allows the
limit of the approximate solutions to pass through the
non-linearity.  To extend to unbounded domains we solve the problem
in balls of radius $\{R_i\}$ (a sequence tending to infinity), and then
invoke a diagonal argument.  Regularity is established
through an inductive argument relying on energy methods.
\end{proof}

Next, we will state a Corollary that describes the action of the filter and will be used many times in the following two sections.

\begin{Corollary}\label{VCHE:filterbounds}
\begin{align}
\|\partial^p_t\nabla^m u\|_2^2 &+ 2\alpha^2 \|\partial^p_t\nabla^{m+1} u\|_2^2+\alpha^4\|\partial^p_t\nabla^{m+2} u\|_2^2
= \|\partial^p_t\nabla^m v\|_2^2\notag\\
\|\partial^p_t\nabla^mu\|_n^2&+ \|\partial^p_t\nabla^{m+1}u\|_n^2 \leq C\|\partial^p_t\nabla^mv\|^2_2\notag\\
\|\partial^p_t \nabla^{m}u(t)\|_n^2 &+\nu\int_0^t\|\partial_t^p \nabla^{m+1}u(s)\|_n^2\,ds \leq C(\|v_0\|_{H^M_0})\notag
\end{align}
for all $m+2k\leq M$, where $C$ is a constant which depends only on $\alpha$,$n$, $m$, and $k$ (in the last bound the constant depends also on $\|v_0\|_{H^M_0}$).
\end{Corollary}
\begin{proof}
This is an application of the Gagliardo-Nirenberg-Sobolev inequality to the bounds in the previous theorem.  Differentiating
the filter relation shows
\begin{equation}
\partial^p_t \nabla^m u-\alpha^2 \partial^p_t \nabla^m \triangle u =\partial^p_t \nabla^m v\notag
\end{equation}
Squaring this relation then integrating by parts gives
\begin{equation}
\|\partial^p_t\nabla^m u\|_2^2 + 2\alpha^2 \|\partial^p_t\nabla^{m+1} u\|_2^2+\alpha^4\|\partial^p_t\nabla^{m+2} u\|_2^2
= \|\partial^p_t\nabla^m v\|_2^2\notag
\end{equation}
This is the first bound in the corollary.  Applying the Gagliardo-Nirenberg-Sobolev inequality to $\|u\|_n$ and using the previous equality shows
\begin{align}
\|\partial^p_t\nabla^mu\|_n^2 \leq C\|\partial^p_t\nabla^mv\|^2_2\notag\\
\|\partial^p_t\nabla^{m+1}u\|_n^2 \leq C\|\partial^p_t\nabla^mv\|^2_2\notag
\end{align}
This is the second set of bounds.  Combining this with the regularity bounds in the theorem give the final set of bounds.
\end{proof}

\section{Large Time Behavior of the VCHE: Non-Uniform Decay}
In bounded domains it is easy to see that the energy of a solution
decays exponentially using the Poincar\'{e} inequality
\begin{equation}
\|u\|_2^2 \leq C(\Omega)\|\nabla u\|_2^2\notag
\end{equation}
Indeed, start with the energy estimate (\ref{VCHE:energy}) and apply
the Poincar\'{e} inequality to find
\begin{equation}
\frac{1}{2}\frac{d}{d t}(<u,u>+\alpha^2<\nabla u, \nabla
u>)+C(\Omega)\nu(<u,u>+\alpha^2<\nabla u,\nabla u>)\leq 0\notag
\end{equation}
This differential inequality implies
\begin{equation}
<u,u>+\alpha^2<\nabla u, \nabla u>\leq C(\|v_0\|_2)e^{-C(\Omega,\nu)t}\notag
\end{equation}
The situation in unbounded domains is more delicate. If the initial data is assumed only in $L^2$ then the solution decays to zero but we are unable to determine the rate without more information, the precise statements of this idea are contained in Theorems \ref{L2deacy:th} and \ref{nonuniform:decaythrm}.

We will follow \cite{MR1432588} to show that the solutions in
the whole space, constructed in Theorem \ref{VCHE:existence}, approach zero as time
becomes large.  The idea is to split the solution into low and high
frequency parts using a cut-off function and generalized energy inequalities to show that both the high and low frequency terms approach zero.  The idea of splitting into low and high
frequency was first used in \cite{MR767409}.

\begin{Lemma}
Solutions of the VCHE constructed in Theorem \ref{VCHE:existence} with $\Omega=\mathbb{R}^n$ satisfy the following generalized energy inequalities.  Let $E\in C^1([0,\infty))$ and $\psi\in C^1([0,\infty);C^1\cap L^2(\mathbb{R}^n))$, then
\begin{align}\label{VCHE:generalizedenergy1}
E(t)\|\psi(t)\ast v(t)\|_2^2 &= E(s)\|\psi(s)\ast v(s)\|_2^2 +\int_s^t E'(\tau) \|\psi(\tau)\ast v(\tau)\|_2^2\,d\tau \\
&+2\int_s^tE(\tau)<\psi'(\tau)\ast v(\tau),\psi(\tau)\ast v(\tau)>\, d\tau \notag\\
&-2\int_s^tE(\tau)\|\nabla\psi(\tau)\ast v(\tau)\|_2^2\, d\tau \notag\\
&-2\int_s^tE(\tau)<u\cdot\nabla v,\psi(\tau)\ast\psi(\tau)\ast v(\tau)>\, d\tau \notag\\
&-2\int_s^tE(\tau)<v\cdot \nabla u^T,\psi(\tau)\ast\psi(\tau)\ast v(\tau)>\, d\tau \notag
\end{align}
For $E\in C^1([0,\infty))$ and $\tilde{\psi}\in C^1(0,\infty;L^\infty(\mathbb{R}^n))$ we have
\begin{align}\label{VCHE:generalizedenergy2}
E(t)\|\tilde{\psi}(t)\hat{v}(t)\|_2^2 &= E(s)\|\tilde{\psi}(s)\hat{v}(s)\|_2^2 +\int_s^t E'(\tau) \|\tilde{\psi}(\tau)\hat{v}(\tau)\|_2^2\,d\tau\\
&+2\int_s^tE(\tau)<\tilde{\psi}'(\tau)\hat{v}(\tau),\tilde{\psi}(\tau)\hat{v}(\tau)>\, d\tau \notag\\
&-2\int_s^tE(\tau)\|\xi\tilde{\psi}(\tau)\hat{v}(\tau)\|_2^2\, d\tau \notag\\
&-2\int_s^tE(\tau)<\mathcal{F}(u\cdot\nabla v),\tilde{\psi}^2(\tau)\hat{v}(\tau)>\, d\tau \notag\\
&-2\int_s^tE(\tau)<\mathcal{F}(v\cdot \nabla u^T),\tilde{\psi}^2(\tau)\hat{v}(\tau)>\, d\tau \notag
\end{align}
\end{Lemma}
\begin{proof}
The proof of the first inequality is accomplished by multiplying the VCHE by $E(t)\psi\ast\psi\ast v$ then integrating by
parts and in time.  The second inequality is obtained by first taking the Fourier Transform of the VCHE, then multiplying by $\tilde{\psi}^2\hat{v}$ and integrating.
\end{proof}

\begin{Theorem}\label{L2deacy:th}
Let $v$ be the solution of the VCHE (\ref{VCHE:PDE}) constructed in Theorem \ref{VCHE:existence} with $\Omega=\mathbb{R}^n$
and $v_0\in L^2_\sigma(\mathbb{R}^n)$, then
\begin{equation}
\lim_{t\rightarrow 0} \|v(t)\|_2 = 0
\end{equation}
\end{Theorem}
\begin{proof}
We work in frequency space . We split the energy into low and high
frequency parts
\begin{equation}
\| \hat{v}\|_2 \leq \|\phi\hat{v}\|_2 + \|(1-\phi)\hat{v}\|_2
\end{equation}
with $\phi=e^{-|\xi|^2}$ will be chosen below.
To estimate the low frequency part of the energy,
begin with the generalized energy estimate (\ref{VCHE:generalizedenergy1}).  Temporarily fix $t$ then choose $E=1$ (the constant function) and
\begin{equation}
\psi(\tau)=\mathcal{F}^{-1}\left[e^{-|\xi|^2(t+1-\tau)}\right]\notag
\end{equation}
Note that $\psi$ and $\mathcal{F}(\psi)$ are rapidly decreasing functions for $\tau<t+1$.  The relation $\hat{\psi}'=|\xi|^2\hat{\psi}$ shows the
third and fourth terms in (\ref{VCHE:generalizedenergy1}) add to zero.
Note $\phi=e^{-|\xi|^2}=\psi(t)$ and apply the Plancherel Theorem to see
\begin{align}\label{VCHE:nonuniformdecayestimate1}
\|\phi\hat{v}(t)\|_2^2&\leq \|e^{|\xi|^2(t-s)}\phi\hat{v}(s)\|_2^2\notag\\
&+2\int_s^t|<\check{\phi}^2\ast (u\cdot\nabla v-v\cdot\nabla u^T),e^{2\triangle(t-\tau)}v(\tau)>|\, d\tau
\end{align}
With H\"{o}lder inequality, Young's inequality, and the Gagliardo-Nirenberg-Sobolev inequality we bound
\begin{align}
|<\check{\phi}^2\ast u\cdot\nabla v,e^{2\triangle(t-\tau)}v(\tau)>|
&\leq \|\check{\phi}^2\ast u\cdot\nabla v\|_2\|e^{2\triangle(t-\tau)}v(\tau)\|_2\notag\\
&\leq C\|\check{\phi}^2\|_\frac{2n}{(n+2)}\|u\|_\frac{2n}{(n-2)}\|\nabla v\|_2\|v\|_2\notag\\
&\leq C(\phi)\|v\|_2\|\nabla u\|_2\|\nabla v\|_2\notag
\end{align}
Similarly,
\begin{align}
|<\check{\phi}^2\ast v\cdot\nabla u^T,e^{2\triangle(t-\tau)}v(\tau)>|
&\leq \|\check{\phi}^2\ast v\cdot\nabla u^T\|_2\|e^{2\triangle(t-\tau)}v(\tau)\|_2\notag\\
&\leq C\|\check{\phi}^2\|_\frac{2n}{(n+2)}\|v\|_\frac{2n}{(n-2)}\|\nabla u\|_2\|v\|_2\notag\\
&\leq C(\phi)\|v\|_2\|\nabla u\|_2\|\nabla v\|_2\notag
\end{align}
Using the triangle inequality, H\"{o}lder's inequality, and (\ref{VCHE:regularityestimate}) in (\ref{VCHE:nonuniformdecayestimate1}) yields
\begin{align}
\|\phi\hat{v}(t)\|_2^2&\leq \|e^{|\xi|^2(t-s)}\phi\hat{v}(s)\|_2^2\notag\\
&+2C(\phi)\|v_0\|_2\left(\int_s^t \|\nabla u\|^2_2\, d\tau\right)^{1/2}\left(\int_s^t\|\nabla v\|^2_2\, d\tau\right)^{1/2}\notag
\end{align}
As the first term on the RHS tends to zero, applying the limit $t\rightarrow \infty$ yields
\begin{align}
\limsup_{t\rightarrow\infty}\|\phi\hat{v}(t)\|_2^2\leq
2C(\phi)\|v_0\|_2\left(\int_s^\infty \|\nabla u\|^2_2\, d\tau\right)^{1/2}\left(\int_s^\infty\|\nabla v\|^2_2\, d\tau\right)^{1/2}\notag
\end{align}
The bounds (\ref{VCHE:aprioriestimate}) and (\ref{VCHE:regularityestimate}) show $\|\nabla u\|_2^2$ and
$\|\nabla v\|_2^2$ are integrable on the positive real line, letting $s\rightarrow\infty$ leaves
\begin{equation}\label{VCHE:lowfrequencydecay}
\limsup_{t\rightarrow\infty}\|\phi\hat{v}(t)\|_2^2\rightarrow 0
\end{equation}

To estimate the high frequency start with the generalized energy inequality (\ref{VCHE:generalizedenergy2}) and chose
$\tilde{\psi}=1-e^{-|\xi|^2}=1-\phi$.  Let $B_G(t)=\{\xi:|\xi|\leq G(t)\}$ where $G(t)$ will be selected later and use $<u\cdot \nabla v, v>=0$ to replace $\tilde{\psi}^2$ by $1- \tilde{\psi}^2$ in the 5th term on the RHS of (\ref{VCHE:generalizedenergy2}).
\begin{align}\label{VCHE:nonuniformdecayestimate2}
E(t)\|\tilde{\psi}\hat{v}(t)\|_2^2 &\leq E(s)\|\tilde{\psi}\hat{v}(s)\|_2^2
+\int_s^t E'(\tau) \int_{B_G(\tau)}|\tilde{\psi}\hat{v}(\tau)|^2\,d\xi\,d\tau\\
&+\int_s^t (E'(\tau)-2E(\tau)G^2(\tau)) \int_{B_G^C(\tau)}|\tilde{\psi}\hat{v}(\tau)|^2\,d\xi\,d\tau\notag\\
&+2\int_s^tE(\tau)|<\mathcal{F}(u\cdot\nabla v+v\cdot \nabla u^T),(1-\tilde{\psi}^2(\tau))\hat{v}(\tau)>|\, d\tau \notag\\
&+2\int_s^tE(\tau)|<\mathcal{F}(v\cdot \nabla u^T),\hat{v}(\tau)>|\, d\tau \notag
\end{align}
We remark both $(1-\tilde{\psi}^2)$ and
$\phi=\mathcal{F}^{-1}(1-\tilde{\psi}^2)$ are rapidly decreasing
functions.  Using again H\"{o}lder's inequality and the Plancherel
theorem, then Young's inequality and the Gagliardo-Nirenberg-Sobolev
inequality allows
\begin{align}
|<\mathcal{F}(u\cdot\nabla v+v\cdot \nabla u^T)&,(1-\tilde{\psi}^2(\tau))\hat{v}(\tau)>| \notag\\
&=|<(1-\tilde{\psi}^2(\tau))\mathcal{F}(u\cdot\nabla v+v\cdot \nabla u^T),\hat{v}(\tau)>| \notag\\
&\leq\|\mathcal{F}^{-1}(1-\tilde{\psi}^2(\tau))\ast(u\cdot\nabla v+v\cdot \nabla u^T)\|_2\|v\|_2\notag\\
&\leq C\|1-\tilde{\psi}\|_\frac{2n}{n+2} (\|u\|_\frac{2n}{n-2}\|\nabla v\|_2+\|v\|_\frac{2n}{n-2}\|\nabla u\|_2)\|v\|_2\notag \\
&\leq C(\phi)\|v\|_2\|\nabla u\|_2\|\nabla v\|_2 \notag
\end{align}
Similarly use H\"{o}lder's inequality with the Plancherel theorem,
then the Gagliardo-Nirenberg-Sobolev inequality, and Corollary \ref{VCHE:filterbounds}
to bound
\begin{align}
|<\mathcal{F}(v\cdot \nabla u^T),\hat{v}(\tau)>|&\leq \|v\cdot\nabla u^T\|_2\|v\|_2\notag\\
&\leq C\|v\|_\frac{2n}{n-2}\|\nabla u\|_n\|v\|_2\notag\\
&\leq C\|v\|_2\|\nabla v\|^2_2\notag
\end{align}
Choosing $E(t)=(1+t)^\beta$ and $G^2(t)=\beta/2(1+t)$ in (\ref{VCHE:nonuniformdecayestimate2}), so that
$E'-2EG^2=0$, and taking $\beta>0$ sufficiently large,  leaves
\begin{align}
\|(1-\phi)\hat{v}(t)\|_2^2 &\leq \frac{(1+s)^\beta}{(1+t)^\beta}\|(1-\phi)\hat{v}(s)\|_2^2\notag\\
&+\int_s^t \frac{\beta(1+\tau)^{\beta-1}}{(1+t)^\beta} \int_{B_G(\tau)}|(1-\phi)\hat{v}(\tau)|^2\,d\xi\,d\tau\notag\\
&+C\|v_0\|_2\int_s^t\frac{(1+\tau)^\beta}{(1+t)^\beta} \|\nabla v\|_2(\|\nabla v\|_2+\|\nabla u\|_2)\, d\tau \notag
\end{align}
For $\xi\in B_G(t)$ and $t$ sufficiently large, $\tilde{\psi}=|1-\phi|\leq |\xi|^2$.  Therefore $|1-\phi|^2\leq\beta^2/4(1+t)^2$
and the second term on the RHS can be bounded as
\begin{align}
\int_s^t\frac{\beta^3(1+\tau)^{-3}}{4}\int_{A(\tau)}|\hat{v}(\tau)|^2\,d\xi\,d\tau\notag
&\leq \int_s^t\frac{\beta^3(1+\tau)^{-3}}{4}\|v(\tau)\|_2^2\,d\tau \notag\\
&\leq \frac{\beta^3}{4}\|v_0\|_2^2\int_s^t(1+\tau)^{-3}\,d\tau\notag\\
&\leq \frac{\beta^3}{8}\|v_0\|_2^2(1+s)^{-2}\notag
\end{align}
Letting $t\rightarrow\infty$ shows
\begin{align}
\limsup_{t\rightarrow\infty}\|(1-\phi)\hat{v}(t)\|_2^2 &\leq \frac{\alpha^3}{8}\|v_0\|_2^2(1+s)^{-2}\\
&+C\|v_0\|_2(\int_s^\infty\|\nabla v\|^2_2\,d\tau + \int_s^\infty\|\nabla u\|^2_2\, d\tau) \notag
\end{align}
The bounds (\ref{VCHE:aprioriestimate}) and (\ref{VCHE:regularityestimate}) again show $\|\nabla v\|^2_2$ and $\|\nabla u\|^2_2$ are integrable on the real line.  Letting $s\rightarrow\infty$ proves
\begin{equation}
\limsup_{t\rightarrow\infty}\|(1-\phi)\hat{v}(t)\|_2^2 =0 \notag
\end{equation}
Combining this with (\ref{VCHE:lowfrequencydecay}) and the Plancherel theorem completes this proof.
\end{proof}
\begin{Corollary}
Let $v$ be the solution of the VCHE (\ref{VCHE:PDE}) constructed in Theorem \ref{VCHE:existence} with $\Omega=\mathbb{R}^n$ corresponding to
$v_0\in H^1_0(\mathbb{R}^n)$.  Then
\begin{equation}
\lim_{t\rightarrow \infty}\frac{1}{t}\int_0^t\|v(\tau)\|_2\,d\tau
=0\notag
\end{equation}
\end{Corollary}
\begin{proof}
Given an $\epsilon>0$ we can choose a large $s$ such that $\|v\|_2\leq \epsilon$ for $\tau>s$, this follows directly from the previous theorem.  Then
\begin{align}
\frac{1}{t}\int_0^t\|v(\tau)\|_2\, d\tau &= \frac{1}{t}\int_0^s\|v(\tau)\|_2\, d\tau +\frac{1}{t}\int_s^t\int \|v(\tau)\|_2 \, d\tau \notag\\
&\leq \frac{1}{t}\int_0^s\|v(\tau)\|_2\, d\tau +\epsilon\frac{t-s}{t}
\end{align}
Note that $\epsilon$ was chosen arbitrarily and let $t\rightarrow \infty$ to finish the proof.
\end{proof}

We have shown that the energy of a solution to the VCHE will tend to zero as time becomes large, now we will provide a
counter example to show that there is no uniform rate of decay based only on the initial energy of the system.
This is analogous to a result proved in \cite{MR837929}.  The idea is to take a family of initial data with a parameter
$\epsilon$ that have constant $L^2$ norm, but norms of higher derivatives of the initial data can be taken arbitrarily small by picking
$\epsilon$ sufficiently small.  It is then possible to bound the higher derivative norms of the solution arbitrarily small by taking $\epsilon$ small..
Combining this with the energy relation (\ref{VCHE:energy}) allows us to place a lower bound on the energy of the solution
which depends on $\epsilon$.  By choosing $\epsilon$ small we can guarantee that a solution will remain away from zero
for any finite amount of time.

\begin{Theorem}\label{nonuniform:decaythrm}
Let $v$ be the solution of the VCHE (\ref{VCHE:PDE}) constructed in Theorem \ref{VCHE:existence} with $\Omega=\mathbb{R}^n$ and $v_0\in L^2_\sigma(\mathbb{R}^n)$.  There exists no function $G(t,\beta):\mathbb{R}^+\times\mathbb{R}^+\rightarrow \mathbb{R}^+$ with the following two properties: 
\begin{align}
\|v\|_2\leq G(t,\|v_0\|_2) \notag\\
\lim_{t\rightarrow\infty}G(t,\beta)= 0 \ \ \ \ \forall\beta
\end{align}
\end{Theorem}
\begin{proof}
Fix $u_0(x)$ to be any smooth function of compact support and write $u^\epsilon_0(x)=\epsilon^{n/2}u_0(\epsilon x)$.
Let $v_0^\epsilon=u^\epsilon_0-\alpha^2\triangle u_0^\epsilon$ and $v^\epsilon$ the solution of the VCHE given by Theorem \ref{VCHE:existence} corresponding to the initial data $v_0$.  Note
$\|u^\epsilon_0\|_2=\|u_0\|_2$ and $\|\nabla^m u^\epsilon_0\|_2=\epsilon^m\|\nabla u_0\|_2$ for all $\epsilon > 0$.  Also,
\begin{align}
\|v_0^\epsilon\|_2^2&=\|u^\epsilon_0\|_2^2+\alpha^2\|\nabla u^\epsilon_0\|_2^2+\alpha^4\|\triangle u^\epsilon_0\|_2^2\\
&=\|u_0\|_2^2+\alpha^2\epsilon^2\|\nabla u_0\|_2^2+\alpha^4\epsilon^4\|\triangle u_0\|_2^2\notag
\end{align}
and
\begin{align}
\|\nabla v_0^\epsilon\|_2^2&=\|\nabla u^\epsilon_0\|_2^2+\alpha^2\|\triangle u^\epsilon_0\|_2^2+\alpha^4\|\nabla\triangle u^\epsilon_0\|_2^2\\
&=\epsilon^2\|\nabla u_0\|_2^2+\alpha^2\epsilon^4\|\triangle u_0\|_2^2+\alpha^4\epsilon^6\|\nabla\triangle u_0\|_2^2\notag
\end{align}
From the two previous inequalities and Corollary \ref{VCHE:filterbounds} we obtain a constant $C=C(\|u_0\|_{H^3_0})$, such that for all $\epsilon >0$
\begin{align}
\|v^\epsilon\|_2^2 &\leq C\label{Gtheorem:eq2}\\
\|\nabla v^\epsilon\|_2^2 &\leq C\epsilon^2 \notag
\end{align}

Multiply the VCHE (\ref{VCHE:PDE}) by $\triangle v^\epsilon$, then
integrating by parts yields
\begin{equation}\label{Gtheorem:eq3}
\frac{1}{2}\frac{d}{dt}\|\nabla v^\epsilon\|^2_2 +\nu \|\triangle^2
v^\epsilon\|^2_2 = <u^\epsilon\cdot\nabla v^\epsilon,\triangle v^\epsilon> + <\triangle
v^\epsilon\cdot\nabla u^\epsilon,v^\epsilon>\notag
\end{equation}
Use the relation $<u^\epsilon,\nabla v^\epsilon,v^\epsilon>=0$, the H\"{o}lder inequality, Sobolev inequality, and then the Cauchy Inequality to see
\begin{align}
|<u^\epsilon\cdot\nabla v^\epsilon,\triangle v^\epsilon>|&= |(-1)<(\nabla u^\epsilon)\cdot\nabla v^\epsilon,\nabla v^\epsilon>| \notag\\
&\leq C\|\nabla u^\epsilon\|_n\|\nabla v^\epsilon\|_2\|\nabla v^\epsilon\|_\frac{2n}{n-2} \notag\\
&\leq \frac{\nu}{4}\|\triangle v^\epsilon\|_2^2 + C \|\nabla u^\epsilon\|^2_n\|\nabla v^\epsilon\|^2_2\notag
\end{align}
Similarly,
\begin{align}
|<\triangle v^\epsilon\cdot\nabla u^\epsilon,v^\epsilon>|&\leq C\|\triangle v^\epsilon\|_2\|\nabla u^\epsilon\|_n \|v^\epsilon\|_\frac{2n}{n-2} \notag\\
&\leq \frac{\nu}{4}\|\triangle v^\epsilon\|_2^2 +C\|\nabla v^\epsilon\|_2^2\|\nabla u^\epsilon\|_n^2 \notag
\end{align}
Applied to (\ref{Gtheorem:eq3}):
\begin{equation}\label{VCHE:CEHELMREL}
\frac{1}{2}\frac{d}{dt}\|\nabla v^\epsilon\|^2_2 +\frac{\nu}{2} \|\triangle
v^\epsilon\|^2_2 \leq C\|\nabla v^\epsilon\|^2_2\|\nabla u^\epsilon\|^2_n
\end{equation}
By (\ref{VCHE:aprioriestimate}) and Corollary \ref{VCHE:filterbounds},
\begin{equation}
\int_0^\infty \|\nabla u^\epsilon\|^2_n \, dt \leq \|u^\epsilon_0\|_2^2+\|\nabla u^\epsilon_0\|_2^2\leq \|v^\epsilon_0\|_2^2\notag
\end{equation}
This bound, combined with (\ref{Gtheorem:eq2}) and (\ref{VCHE:CEHELMREL}) yields
\begin{equation}
\|\nabla v^\epsilon\|_2^2 \leq \|\nabla v^\epsilon_0\|^2_2e^{C\|v^\epsilon_0\|_2^2}\leq C\epsilon^2 e^{C\epsilon^2}\notag
\end{equation}
Again, apply Corollary \ref{VCHE:filterbounds}
\begin{equation}
\|\nabla u^\epsilon\|_2^2 +\alpha^2 \|\triangle u^\epsilon\|_2^2\leq\|\nabla v^\epsilon\|_2^2 \leq \|\nabla v^\epsilon_0\|^2_2e^{C\|v^\epsilon_0\|_2^2}\notag
\end{equation}
This together with the energy estimate (\ref{VCHE:energy}) implies
\begin{equation}
\frac{1}{2}\frac{d}{dt}(\|u^\epsilon\|_2^2+\alpha^2\|\nabla u^\epsilon\|_2^2)\geq -C\epsilon^2 \notag
\end{equation}
or,
\begin{align}\label{VCHE:CEHELMREL2}
\|u^\epsilon\|_2^2+\alpha^2\|\nabla u^\epsilon\|_2^2 &\geq \|u^\epsilon_0\|_2^2+\alpha^2\|\nabla u^\epsilon_0\|_2^2 -C\epsilon^2 t\\
&= \|u_0\|_2^2+\epsilon^2\alpha^2\|\nabla u_0\|_2^2 -C\epsilon^2 t \notag \\
&\geq \|u_0\|_2^2-C\epsilon^2 t \notag
\end{align}
From this we can deduce that there is no function $G(t,\beta,)$,
continuous and approaching zero in $t$ for each fixed $\beta$, such that $\|u\|_2 \leq G(t,\|u_0\|_2)$.  If
there was such a function, then at some $t_0$ it would satisfy the
bound $G(t_0,\|u_0\|_2)\leq \|u_0\|_2/2$.  By choosing $\epsilon$
sufficiently small in (\ref{VCHE:CEHELMREL2}), i.e. $\epsilon^2
<\|u_0\|_2^2/4Ct_0 $, we have found initial data with a solution
which cannot satisfy this estimate.
\end{proof}

\section{Large Time Behavior of the VCHE: Algebraic Decay}
Although there is no uniform rate of decay for solutions with data exclusively in $L^2$, we now show that
there is a uniform rate of decay depending on the $L^2$ and $L^1$ norm of the initial data.  Theorem \ref{VCHE:decay} contains the most general decay result in this section.  

There is a relation between the shape of the Fourier Transform of the initial data near the origin and the decay rate of a solution to a parabolic equation with this data.  By requiring the initial data to be absolutely integrable (in $L^1$) we are in turn requiring the Fourier Transform of the initial data to be bounded.
Using the Fourier Splitting Method it will also be shown that solutions in
the whole space decay algebraically  in $H^M$ as $t\rightarrow \infty$ for
initial data in $L^1\cap H^M$, $M\geq 0$.  The decay obtained is the same as for the linear part (the heat equation).   Note
that the initial conditions can be weakened to require only that $v_0\in X$ where $X=\{v_0|v_0(t)\leq C(1+t)^{-\beta}\}$ where
$v_0(t)$ is the solution of the heat equation with initial data $v_0$.  The decay rate will depend on the relation between $\beta$
and the number of dimensions.  For similar results corresponding to the Navier-Stokes equations see \cite{MR1396285}, \cite{MR881519}, and \cite{MR1993060}.

The Fourier Splitting Method was originally applied to parabolic conservation laws in \cite{FSPCL}, and later applied to the NSE in
\cite{MR775190}.  In \cite{MR837929} the decay rate was made sharp
in dimension $n>2$ through a bootstrap method and logarithmic decay
was shown for $n=2$.  In \cite{MR1312702} the decay rate for $n=2$
was made sharp through a bootstrap argument involving the Gronwall
inequality.  In this section we combine ideas from all of these
papers in a slightly different way which allows us to prove the
optimal energy decay rate in dimensions $n\geq 2$ without appealing to a
bootstrap argument although we still use a bootstrap argument to obtain decay rates for higher derivatives.  This same argument is also applicable to the
NSE.

The first goal of this section is to obtain a decay rate for the
filtered velocity $u$, which is accomplished by applying the Fourier
Splitting Method to the natural energy relation (\ref{VCHE:energy}).
This decay rate is then used with an inductive argument to obtain decay
rates for the unfiltered velocity $v$ and all of its derivatives.
We start by finding estimates on $\|\hat{v}\|_\infty$.

\begin{Lemma}\label{VCHE:vfrequencybound}
Let $v$ be the solution of the VCHE (\ref{VCHE:PDE}) constructed in Theorem \ref{VCHE:existence} with $\Omega=\mathbb{R}^n$, corresponding to
$v_0\in L^2_\sigma\cap L^1(\mathbb{R}^n)$.
Then,
\begin{equation}
|\mathcal{F}(v)|\leq C\left[1+\left(\int_0^t\|u(s)\|^2_2\right)^{1/2}\left(\int_0^t\|\nabla v(s)\|^2_2\right)^{1/2}\right] \notag
\end{equation}
where the constant depends only on the initial data, the dimension of space, and the constants in the VCHE (but not $\alpha$).
\end{Lemma}
\begin{proof}
Use the identity
\begin{equation} \label{cheapidentity}
\sum_i\nabla(u_iv_i)=\sum_i u_i\nabla v_i +\sum_i v_i\nabla u_i
\end{equation}
and write the Fourier transform of the
solution $\mathcal{F}(v)$ as
\begin{equation}\label{mildFT}
\mathcal{F}(v)=e^{-\nu t|\xi|^2}\mathcal{F}(v_0)
+\int^t_0e^{-\nu(t-s)|\xi|^2}\Psi(\xi,s)\,ds \notag
\end{equation}
where
\begin{equation}\label{mildFT2}
\Psi(\xi,t)=\xi\cdot\mathcal{F}(\pi+\sum_i u_iv_i)- \mathcal{F}(u\cdot\nabla v - u\cdot\nabla v^T)
\end{equation}
We would like first to bound $\Psi$, in that direction we have the following estimate which relies on the bound $\|\mathcal{F}(\phi)\|_\infty\leq \|\phi\|_1$ and Young's inequality
\begin{align}
|\mathcal{F}(u\cdot\nabla v - u\cdot\nabla v^T)|\leq C\|u\|_2\|\nabla v\|_2 \notag
\end{align}
Also, taking the divergence of the VCHE (\ref{VCHE:PDE}) shows
\begin{equation}
\triangle (\pi+\sum_i u_iv_i) = \div (u\nabla v - u\nabla v^T) \notag
\end{equation}
Using the estimate immediately above and the Fourier transform leaves
\begin{equation}
|\xi \mathcal{F}(\pi+\sum_i u_iv_i)|\leq C\|u\|_2\|\nabla v\|_2\notag
\end{equation}
Now we can bound the integrand
\begin{equation}
|\Psi(\xi,t)|\leq C\|u\|_2\|\nabla v\|_2\notag
\end{equation}
Now take the supremum over $\xi$ of (\ref{mildFT}) and apply the Cauchy-Schwartz Inequality:
\begin{equation}\label{vfreq:cor1}
|\mathcal{F}(v)|\leq|\mathcal{F}(v_0)| + C\left(\int^t_0\|u(s)\|_2^2\,ds\right)^{1/2}\left(\int^t_0\|\nabla v(s)\|_2^2\,ds\right)^{1/2} \notag
\end{equation}
The bound $|\mathcal{F}(v_0)|_\infty\leq \|v_0\|_1$ finishes the proof.
\end{proof}

\begin{Theorem}\label{VCHE:decaytheorem1}
Let $v$ be the solution of the VCHE (\ref{VCHE:PDE}) constructed in Theorem \ref{VCHE:existence} with $\Omega=\mathbb{R}^n$, corresponding to
$v_0\in L^2_\sigma\cap L^1(\mathbb{R}^n)$.  The solution satisfies the ``energy'' decay rate
\begin{equation}
\int_{\mathbb{R}^n} v\cdot u\,dx = \|u\|^2_2+\alpha^2\|\nabla
u\|^2_2 \leq C (t+1)^{-n/2} \notag
\end{equation}
where the constant depends only on the initial data, the dimension
of space, and the constants in the VCHE (but not $\alpha$).
\end{Theorem}
\begin{proof}
The previous lemma, with the bound (\ref{VCHE:regularityestimate}), yields
\begin{equation}\label{bound:Fv1}
|\hat{v}|^2\leq C\left[1+\int^t_0\|u(s)\|_2^2\,ds\right]
\end{equation}

Now we begin work with the energy estimate (\ref{VCHE:energy}).
Using the Plancherel Theorem we rewrite it as
\begin{equation}
\frac{d}{dt}\int_{\mathbb{R}^n} (1+\alpha^2 |\xi|^2)\hat{u}^2\,d\xi
+2\nu\int_{\mathbb{R}^n} |\xi|^2(1+\alpha^2
|\xi|^2)\hat{u}^2\,d\xi=0 \notag
\end{equation}
Let $B(\rho)$ be the ball of radius $\rho$ where $\rho^2=f'(t)/(2\nu
f(t))$, and $f$ is a positive, increasing function to be specified
later.  To simplify our equations we write $E^2=\hat{u}\cdot\hat{v}=(1+\alpha^2 |\xi|^2)\hat{u}^2$.
Then,
\begin{equation}
\frac{d}{dt}\int_{\mathbb{R}^n} E^2\,d\xi
+2\nu\rho^2\int_{B^C(\rho)}E^2\,d\xi\leq 0 \notag
\end{equation}
or
\begin{equation}\label{E:E1}
\frac{d}{dt}\int_{\mathbb{R}^n} E^2\,d\xi +2\nu\rho^2\int_{\mathbb{R}^n}E^2\,d\xi\leq 2\nu\rho^2\int_{B(\rho)}E^2\,d\xi
\end{equation}

Recall the relation between $u$ and $v$, that is
$v=u-\alpha^2\triangle u$ which has Fourier Transform $\hat{u}=\hat{v}/(1+\alpha^2|\xi|^2)$.  Combining this with (\ref{bound:Fv1}) we see
\begin{align}
\|E\|_\infty^2 &= \frac{\|\hat{v}\cdot\hat{v}\|_\infty}{(1+\alpha^2|\xi|^2)} \notag\\
&\leq C\left[1+\int^t_0\|u(s)\|_2^2\,ds\right] \notag
\end{align}
With this bound we can estimate the integral on the right hand side of (\ref{E:E1}).
\begin{equation}\label{E:E3}
\frac{d}{dt}\int_{\mathbb{R}^n} E^2\,d\xi
+2\nu\rho^2\int_{\mathbb{R}^n}E^2\,d\xi \leq C\rho^{2+n} \left[1+\int^t_0\|u(s)\|_2^2\,ds\right] \notag
\end{equation}
We now have a differential inequality which can be solved using the
integrating factor $f$ to find
\begin{equation}
\frac{d}{dt}\left(f\int_{\mathbb{R}^n} E^2\,d\xi\right)\leq
C f'\left(\frac{f'}{f}\right)^{n/2}\left[1+\int^t_0\|u(s)\|_2^2\,ds\right]\notag
\end{equation}
Choose $f=(1+t)^{n/2+1}$ so that $f'/f=(n/2+1)/(1+t)$and integrate in time from $0$ to $r$.
\begin{equation}
(1+r)^{n/2+1}\int_{\mathbb{R}^n} E^2(\xi,r)\,d\xi\leq
\int_{\mathbb{R}^n}
E^2(\xi,0)\,d\xi+C\int_0^r(1+\int_0^t\|u(s)\|_2^2\,ds)\,dt\notag
\end{equation}

Note $\|u\|_2^2\leq \int_{\mathbb{R}^n} E^2\,d\xi$, then using the Tonelli theorem we can bound the integral on the RHS as
\begin{align}
\int_0^r(1+\int_0^t\|u(s)\|_2^2\,ds)\,dt&\leq  \int_0^r(1+\int_0^r\|u(s)\|_2^2\,ds)\,dt\notag\\
&\leq r(1+\int_0^r\int_{\mathbb{R}^n} E^2(\xi,s)\,d\xi\,ds)\notag
\end{align}
which leaves
\begin{align}
(1+r)^{n/2+1}\int_{\mathbb{R}^n} E^2(\xi,r)\,d\xi &\leq \int_{\mathbb{R}^n} E^2(\xi,0)\,d\xi+Cr(1+\int_0^r\|u\|_2^2\,ds)\notag\\
&\leq C(1+r)+Cr\int_0^r\int_{\mathbb{R}^n} E^2(\xi,s)\,d\xi\,ds\notag
\end{align}
This is of the form 
\begin{equation}
\phi\leq C(1+r)+C(r)\int_0^r \phi(s)(1+s)^{-n/2+1}\,ds\notag
\end{equation}
with $\phi=(1+r)^{n/2+1}\int_{\mathbb{R}^n} E^2(\xi,r)\,d\xi$.
The Gronwall inequality now shows
\begin{equation}\notag
(1+r)^{n/2+1}\int_{\mathbb{R}^n} E^2(\xi,r)\,d\xi\leq
C(1+r)\exp(Cr\int_0^r(1+s)^{-n/2-1}\,ds)
\end{equation}
For $n\geq 2$ the integral $r\int_0^r(1+s)^{-n/2-1}\,ds$ is bounded independent of $r$.  Applying the Plancherel theorem one more time finishes the proof.
\end{proof}

Next we work out of order and establish the decay rate for the homogeneous $H^1$ norm of $v$ using a similar argument as the previous theorem.
\begin{Theorem}\label{VCHE:decaytheorem2}
Let $v$ be the solution of the VCHE (\ref{VCHE:PDE}) constructed in
Theorem \ref{VCHE:existence} with $\Omega=\mathbb{R}^n$
corresponding to $v_0\in H^1_\sigma\cap L^1(\mathbb{R}^n)$.  The
solution satisfies the decay rate
\begin{equation}
\|\nabla v\|^2_2 \leq C (t+1)^{-1-n/2} \notag
\end{equation}
where the constant depends only on the initial data, the dimension
of space, and the constants in the VCHE ($\nu$, $\alpha$).
\end{Theorem}
\begin{proof}
Multiply the VCHE by $\triangle v$, use the identity (\ref{cheapidentity}), after recalling that $\triangle v$ is divergence free use the H\"{o}lder inequality to obtain
\begin{equation}
\frac{1}{2}\frac{d}{dt}\|\nabla v\|_2^2 +\nu\|\triangle v\|_2^2\leq C\|u\|_n\|\nabla v\|_\frac{2n}{n-2}\|\triangle v\|_2\notag
\end{equation}
After using the Sobolev inequality, Corollary \ref{VCHE:filterbounds}, and the previous theorem, this becomes
\begin{equation}
\frac{1}{2}\frac{d}{dt}\|\nabla v\|_2^2 +\nu\|\triangle v\|_2^2\leq C(1+t)^{-n/2}\|\triangle v\|^2_2\notag
\end{equation}
We will now restrict ourselves to $t$ large enough so that $C(1+t)^{-1}<\nu/2$, this implies
\begin{equation}
\frac{d}{dt}\|\nabla v\|_2^2 +\nu\|\triangle v\|_2^2\leq 0 \notag
\end{equation}
The next step is to apply the Fourier Splitting method as in the previous theorem.
Let $B(\rho)$ be the ball of radius $\rho$ where $\rho^2=f'/(\nu f)$ and $f$ is a positive increasing function to be specified later,
using the Plancherel theorem:
\begin{equation}
\frac{d}{dt}\|\xi \hat{v}\|_2^2 +\nu\rho^2\|\xi \hat{v}\|_2^2\leq \nu\rho^4\int_{B(\rho)}|\hat{v}|^2\,\xi\notag
\end{equation}
Lemma \ref{VCHE:vfrequencybound} with Theorem \ref{VCHE:decaytheorem1} imply
\begin{equation}
|\hat{v}|^2\leq C\left[1+\left(\int_0^t (1+s)^{-n/2}\,ds\right)\left(\int_0^t\|\nabla v\|_2^2\,ds\right)\right]\notag\\
\end{equation}
With this bound the previous line becomes
\begin{equation}
\frac{d}{dt}\|\xi \hat{v}\|_2^2 +\nu\rho^2\|\xi \hat{v}\|_2^2
\leq C\nu\rho^{4+n}\left[1+\left(\int_0^t (1+s)^{-n/2}\,ds\right)\left(\int_0^t\|\nabla v\|_2^2\,ds\right)\right]\notag
\end{equation}
Set $f=(1+t)^{n/2+2}$ and use it as an integrating factor
\begin{equation}
\frac{d}{dt}\left((1+t)^{n/2+2}\|\xi \hat{v}\|_2^2\right)\leq
C\left[1+\left(\int_0^t (1+s)^{-n/2}\,ds\right)\left(\int_0^t\|\nabla v\|_2^2\,ds\right)\right]\notag
\end{equation}
Again, as in the previous theorem, integrate in time from $0$ to $r$, then use the Tonelli theorem and the Plancherel theorem
\begin{equation}
(1+r)^{n/2+2}\|\nabla v\|_2^2\leq  C(1+r) \left(\int_o^r \left(\int_0^t (1+s)^{-n/2}\,ds\right)\,dt\right) \int_0^r\|\nabla v(s)\|^2_2\,ds\notag
\end{equation}
The Gronwall inequality now shows
\begin{equation}
(1+r)^{n/2+2}\|\nabla v\|_2^2\leq C(1+r)e^A\notag
\end{equation}
Where
\begin{equation}
A=\left[\left(\int_o^r \left(\int_0^t (1+s)^{-n/2}\,ds\right)\,dt\right)\left(\int_o^r (1+t)^{-n/2-2}\,dt\right)\right]\notag
\end{equation}
Note $A$ is finite, hence,
\begin{equation}
\|\nabla v(r)\|_2^2\leq C(1+r)^{-n/2-1}\notag
\end{equation}
\end{proof}

\begin{Corollary}\label{VCHE:specboundcorr}
Let $v$ be the solution of the VCHE (\ref{VCHE:PDE}) constructed in Theorem \ref{VCHE:existence} with $\Omega=\mathbb{R}^n$, corresponding to
$v_0\in L^2_\sigma\cap L^1(\mathbb{R}^n)$.  Then,
\begin{align}
|\mathcal{F}(v)|\leq C \notag\\
|\mathcal{F}(u)|\leq C  \notag
\end{align}
where the constant depends only on the initial data, the dimension
of space, and the constants in the VCHE.
\end{Corollary}
\begin{proof}
Combine Lemma \ref{VCHE:vfrequencybound} with Theorems
\ref{VCHE:decaytheorem1} and \ref{VCHE:decaytheorem2}.
\end{proof}

\begin{Corollary}\label{VCHE:vdecay}
Let $v$ be the solution of the VCHE (\ref{VCHE:PDE}) constructed in Theorem \ref{VCHE:existence} with $\Omega=\mathbb{R}^n$, corresponding to
$v_0\in L^2_\sigma\cap L^1(\mathbb{R}^n)$.  Then
\begin{equation}
\|v\|^2_2\leq C (t+1)^{-n/2} \notag
\end{equation}
where the constant depends only on the initial data, the dimension
of space, and the constants in the VCHE.
\end{Corollary}
\begin{proof}
In Theorem \ref{VCHE:decaytheorem1} we have shown that
\begin{equation}\label{VCHE:vdecayref1}
\|u\|^2_2+\alpha^2\|\nabla u\|^2_2 \leq C (t+1)^{-n/2}
\end{equation}
Differentiating the Helmholtz equation, then squaring it and integrating shows
\begin{equation}
\|\nabla u\|^2_2+2\alpha^2\|\nabla^2 u\|^2_2+ \alpha^4\|\nabla^3 u\|^2_2 =\|\nabla v\|_2^2  \notag
\end{equation}
Combine this with Theorem \ref{VCHE:decaytheorem2},
\begin{equation}
\|\nabla^2 u\|^2_2 \leq C (t+1)^{-n/2-1}  \notag
\end{equation}
With (\ref{VCHE:vdecayref1}) we see
\begin{equation}
\|v\|_2^2 = \|u\|_2^2+2\alpha^2\|\nabla u\|_2^2 + \alpha^4\|\triangle u\|_2^2 \leq C(t+1)^{-n/2}\notag
\end{equation}
\end{proof}

We now turn our attention to a more general situation involving the
Fourier Splitting Method.  This next theorem will be used in the
remaining decay proofs.

\begin{Theorem}\label{FSM}
Let $\|\nabla^m w(0)\|_2<\infty$.  Given an energy inequality of the form
\begin{equation}
\frac{1}{2}\frac{d}{dt}\|\nabla^m w\|^2_2 + \nu\|\nabla^{m+1}
w\|^2_2 \leq C(1+t)^\gamma \notag
\end{equation}
and the bound
\begin{equation}
|\hat{w}(\xi,t)|\leq C(1 +t)^\beta \notag
\end{equation}
which holds for $|\xi|^2<\frac{d}{\nu(1+t)}$, we can deduce the
asymptotic behavior
\begin{equation}
\|\nabla^{m} w\|^2_2 \leq C\left[(1+t)^{-m-n/2+2\beta} +
(1+t)^{\gamma+1}\right] \notag
\end{equation}
\end{Theorem}
\begin{proof}
We proceed directly with Fourier Splitting.  Apply Plancherel's
Theorem and break up the integral on the LHS.
\begin{equation}
\frac{1}{2}\frac{d}{dt}\|\xi^k \hat{w}\|^2_2 + \nu\rho^2\|\xi^{k}
\hat{w}\|^2_2 \leq \nu\rho^{2k+2}\int_{B(\rho)}\hat{w}^2\,d\xi +
C(1+t)^\gamma \notag
\end{equation}
Choose, for some large $d$,
\begin{equation}
\rho^2=\frac{d}{\nu (1+t)} \notag
\end{equation}
Then, using the assumption for the bound on $\hat{w}$ and performing
the integration on the RHS we have
\begin{equation}
\frac{d}{dt}((1+t)^d\|\xi^m \hat{w}\|^2_2)\leq
C\left[(1+t)^{-m-1+d+2\beta-n/2}+(1+t)^{\gamma+d}\right] \notag
\end{equation}
Integration in time and another application of the Plancherel
theorem  finishes the proof.
\end{proof}

For the first application of the above theorem we will compute the decay rate for all spacial derivatives for solutions of the VCHE.

\begin{Theorem}\label{VCHE:spacedecay}
Let $v$ be the solution of the VCHE (\ref{VCHE:PDE}) constructed in Theorem \ref{VCHE:existence} with $\Omega=\mathbb{R}^n$, corresponding to
$v_0\in H^K_\sigma\cap L^1(\mathbb{R}^n)$.  These solutions satisfy the following decay
rate for all $m\leq K$
\begin{equation}
\|\nabla^m v\|^2_2\leq C (t+1)^{-m-n/2} \notag
\end{equation}
\end{Theorem}
\begin{proof}
The cases $m=0,1$ are Theorems \ref{VCHE:vdecay} and \ref{VCHE:decaytheorem2} respectively.  To prove the remaining cases, we first find an inequality in a form suitable for Theorem \ref{FSM}, then using inductive arguments establish decay.  Having previously
established regularity of solutions, we proceed formally.  Let $M\leq K$ then
multiply the VCHE (\ref{VCHE:PDE}) by $\triangle^{M}v$
and integrate by parts to find
\begin{equation}
\frac{d}{dt}\|\nabla^{M} v\|^2_2+\nu\|\nabla^{M+1} v\|_2^2 \leq I_{M,0} +J_{M,0} \notag
\end{equation}
where
\begin{align}
I_{M,0}=\sum_{m=0}^{M}{M\choose m}
<\nabla^mu\cdot \nabla^{M+1}v,\nabla^{M-m}v> \notag\\
J_{M,0}=\sum_{m=0}^{M-1}{M-1\choose m}
<\nabla^{M+1}v\cdot\nabla^{m+1}u,\nabla^{M-m}v> \notag
\end{align}
Using the Holder inequality, the Sobolev inequality, Corollary \ref{VCHE:filterbounds}, and the Cauchy inequality we find
\begin{align}
I_{M,0}&=C\sum_{m=0}^{M}
\|\nabla^mu\|_n\|\nabla^{M-m}v\|_\frac{2n}{n-2}\|\nabla^{M+1}v\|_2 \notag\\
&\leq  C\| v\|_2^2\|\nabla^{M+1} v\|_2^2 +C\|\nabla v\|_2^2\|\nabla^Mv\|_2^2 \notag\\
&+C\sum_{m=2}^{M} \|\nabla^{m-1}v\|^2_2\|\nabla^{M+1-m}v\|^2_2 +\frac{\nu}{4}\|\nabla^{M+1}v\|^2_2\notag
\end{align}
We treat the other term in a similar way.
\begin{align}
J_{M,0} &\leq C\sum_{m=0}^{M-1} \|\nabla^{M+1}v\|_2\|\nabla^{m+1}u\|_n\|\nabla^{M-m}v\|_\frac{2n}{n-2} \notag\\
&\leq C\| v\|_2^2\|\nabla^{M+1} v\|_2^2+C\|\nabla v\|_2^2\|\nabla^Mv\|_2^2 \notag\\
&+C\sum_{m=2}^{M-1} \|\nabla^{m-1}v\|^2_2\|\nabla^{M+1-m}v\|^2_2 +\frac{\nu}{4}\|\nabla^{M+1}v\|^2_2\notag
\end{align}
Together, this leaves
\begin{align}\label{spacedecayinductionstart}
\frac{d}{dt}\|\nabla^{M} v\|^2_2+\frac{\nu}{2}\|\nabla^{M+1} v\|_2^2 &\leq C\| v\|_2^2\|\nabla^{M+1} v\|_2^2\\
&+C\|\nabla v\|_2^2\|\nabla^Mv\|_2^2 \notag\\
&+C\sum_{m=2}^{M} \|\nabla^{m-1}v\|^2_2\|\nabla^{M+1-m}v\|^2_2 \notag
\end{align}
The remaining part of this proof will proceed by induction where the base case is Theorems \ref{VCHE:decaytheorem2} and \ref{VCHE:vdecay}.  We assume (inductive assumption) that the decay
\begin{equation}
\|\nabla^m v\|^2_2\leq C (t+1)^{-m-n/2} \notag
\end{equation}
holds for all $m< M$ and will show that it holds for $m=M$.  With this inductive assumption (\ref{spacedecayinductionstart}) becomes
\begin{align}\label{spacedecayinductionmid}
\frac{d}{dt}\|\nabla^{M} v\|^2_2+\frac{\nu}{4}\|\nabla^{M+1} v\|_2^2 &\leq C(1+t)^{-n/2}\|\nabla^{M+1} v\|_2^2\\
&+C(1+t)^{-1-n/2}\|\nabla^Mv\|_2^2 +C(1+t)^{-M-n/2}\notag
\end{align}
Consider $t$ large enough so that $C(1+t)^{-n/2} \leq \nu/4$.  Subtracting the first term on the RHS, (\ref{spacedecayinductionmid}) becomes
\begin{align}\label{spacedecayinductionend}
\frac{d}{dt}\|\nabla^{M} v\|^2_2+\frac{\nu}{4}\|\nabla^{M+1} v\|_2^2 &\leq C(1+t)^{-1-n/2}\|\nabla^Mv\|_2^2 \notag\\
&+C(1+t)^{-M-n/2}\notag
\end{align}
The next step is to apply the bound $\|\nabla^Mv\|_2^2\leq C$ (Theorem \ref{VCHE:existence}) with Theorem \ref{FSM} to obtain the decay rate
\begin{equation}
\|\nabla^Mv\|_2^2\leq C(1+t)^{-n/2}\notag
\end{equation}
Continuing with a bootstrap argument, placing this new bound into (\ref{spacedecayinductionmid}) and again using Theorem \ref{FSM} the optimal decay rate is obtained and the proof is complete.
\end{proof}

The next goal is to extend the decay results to time derivatives of the solution.  To begin we will compute a frequency bound for the spacial derivatives
of solutions to the VCHE.  This next lemma will be used inductively with Theorem
\ref{FSM} to compute decay rates for the $L^2$ norm of all time derivatives.

\begin{Lemma}\label{VCHE:vtfrequencybound}
Let $P\geq 1$ and $v$ be the solution of the VCHE (\ref{VCHE:PDE}) constructed in Theorem \ref{VCHE:existence} with $\Omega=\mathbb{R}^n$, corresponding to
$v_0\in H^1_\sigma\cap L^1(\mathbb{R}^n)$.  If
\begin{equation}
\|\partial^p_t \nabla^m v\|^2_2\leq C(1+t)^{-2p-m-n/2} \notag
\end{equation}
for all $p<P$ and $m=0,1$, then
\begin{equation}
|\partial^P_t\hat{v}(\xi)|\leq C(1+t)^{-P} \notag
\end{equation}
for $|\xi|^2 \leq \frac{d}{\nu(1+t)}$.
Here the constant depends only on the initial data, the dimension
of space, and the constants in the VCHE.
\end{Lemma}
\begin{Remark}
Note that the conclusion for $P=0$ is true by Corollary \ref{VCHE:specboundcorr}.
\end{Remark}
\begin{proof}
The chain rule
\begin{equation}
\frac{d}{dt} \int_0^t f(t,s) \,ds = f(t,t)+\int_0^t\frac{\partial
f(t,s)}{\partial t}\,ds \notag
\end{equation}
applied to (\ref{mildFT}) shows
\begin{align}
\partial^P_t \mathcal{F}(v) = (-1)^P|\xi|^{2P}e^{-t|\xi|^2}\mathcal{F}(v_0)&+
\sum_{p=0}^{P-1} (-|\xi|^2)^{P-1-p} \partial^p_t \Psi(\xi,t) \notag\\
&+ \int_0^t(-|\xi|^2)^P e^{-(t-s)|\xi|^2}\Psi(\xi,s)\,ds \notag
\end{align}
We bound $\Psi$ (defined by \ref{mildFT2}) similar to the proof of Lemma \ref{VCHE:vfrequencybound} but
using the assumptions of this Lemma.
\begin{align}
\partial^p_t \Psi(\xi,t)=\partial^p_tA+\partial^p_tB+\partial^p_tC \notag
\end{align}
\begin{align}
|\partial^p_tA|&=|\partial^p_t\sum_j \xi_j\mathcal{F}(u_jv)| \notag\\
&\leq \sum_{l=0}^pC|\xi|\|\partial^l_t v\|_2\|\partial^{p-l}_t v\|_2 \notag\\
&\leq C(1+t)^{-p-n/2-1/2} \notag
\end{align}
\begin{align}
|\partial^p_tB|&=|\partial^p_t\sum_j \mathcal{F}(u_j\nabla v_j^T)| \notag\\
&\leq \sum_{l=0}^pC\|\partial^l_t v\|_2\|\partial^{p-l}_t \nabla v\|_2 \notag\\
&\leq C(1+t)^{-p-n/2-1/2} \notag
\end{align}
\begin{align}
|\partial^p_tC|&=|\partial_t^p\xi\mathcal{F}(\pi)| \notag\\
&\leq |\partial_t^pA|+|\partial_t^pB| \notag\\
&\leq C(1+t)^{-p-n/2-1/2} \notag
\end{align}
The bound $|\hat{v}|\leq C$ (Corollary
\ref{VCHE:specboundcorr}) and $|\xi|<\frac{d}{\sqrt{\nu(1+t)}}$ finish the proof.
\end{proof}

\begin{Theorem}\label{VCHE:decay}
Let $v$ be the solution of the VCHE (\ref{VCHE:PDE}) constructed in Theorem \ref{VCHE:existence} with $\Omega=\mathbb{R}^n$, corresponding to
$v_0\in H^K_\sigma\cap L^1(\mathbb{R}^n)$.  These solutions satisfy the following decay
rate for all $m+2p\leq K$
\begin{equation}
\|\partial^p_t\nabla^m v\|^2_2\leq C (t+1)^{-2p-m-n/2} \notag
\end{equation}
where the constant depends only on the initial data, the dimension
of space, and the constants in the VCHE.
\end{Theorem}
\begin{proof}
This proof follows closely the proof of Theorem \ref{VCHE:spacedecay}, we first find an inequality in a form suitable for Theorem \ref{FSM}, then using inductive arguments we establish decay.  Choose $P$ and $M$ such that $M+2P\leq K$, then apply $\partial_t^P$ to the VCHE (\ref{VCHE:PDE}), multiply by $\partial_t^P\triangle^{M}v$
and integrate by parts to see
\begin{equation}
\frac{d}{dt}\|\partial^{P}_t \nabla^{M} v\|^2_2+\nu\|\partial^{P}_t
\nabla^{M+1} v\|_2^2 \leq I_{M,P} +J_{M,P} \notag
\end{equation}
where
\begin{align}
I_{M,P}&=\sum_{p=0}^{P}\sum_{m=0}^{M}{P\choose p}{M\choose m}
<\partial^p_t\nabla^mu\cdot \partial^{P}_t\nabla^{M+1}v,\nabla^{M-m}\partial^{P-p}_tv> \notag\\
J_{M,P}&=\sum_{p=0}^{P}\sum_{m=0}^{M-1}
{P\choose p}{M-1\choose m}
<\partial^{P}_t\nabla^{M+1}v\cdot\nabla\partial^{p}_t\nabla^{m}u,\partial^{P-p}_t\nabla^{M-m}v> \notag
\end{align}
or, in the case $M=0$,
\begin{equation}
J_{0,P}=\sum_{p=0}^{P}{P\choose p}<\partial^{P}_tv\cdot\nabla\partial^{p}_t u,\partial^{P-p}_tv> \notag
\end{equation}

Use the Holder inequality, the Sobolev inequality, Corollary \ref{VCHE:filterbounds}, and the Cauchy inequality we find, for $M>0$,
\begin{align}
I_{M,P}&=C\sum_{p=0}^{P}\sum_{m=0}^{M}
\|\partial^p_t\nabla^mu\|_n\|\partial^{P-p}_t\nabla^{M-m}v\|_\frac{2n}{n-2}\|\partial^{P}_t\nabla^{M+1}v\|_2 \notag\\
&\leq C\sum_{p=0}^{P}\|\partial_t^p v\|_2^2\|\partial_t^{P-p}\nabla^{M+1} v\|_2^2+C\sum_{p=0}^P\|\partial_t^p\nabla v\|_2^2\|\partial_t^{P-p}\nabla^Mv\|_2^2 \notag\\ &+C\sum_{p=0}^{P}\sum_{m=2}^{M} \|\partial^p_t\nabla^{m-1}v\|^2_2\|\partial^{P-p}_t\nabla^{M+1-m}v\|^2_2 +\frac{\nu}{4}\|\partial^{P}_t\nabla^{M+1}v\|^2_2\notag
\end{align}
Similarly for the second term if $M>0$,
\begin{align}
J_{M,P}&\leq
C\sum_{p=0}^{P}\sum_{m=0}^{M-1} \|\partial^{P}_t\nabla^{M+1}v\|_2\|\partial^p_t\nabla^{m+1}u\|_n\|\partial^{P-p}_t\nabla^{M-m}v\|_\frac{2n}{n-2} \notag\\
&\leq C\sum_{p=0}^{P}\|\partial_t^p v\|_2^2\|\partial_t^{P-p}\nabla^{M+1} v\|_2^2+C\sum_{p=0}^P\|\partial_t^p\nabla v\|_2^2\|\partial_t^{P-p}\nabla^Mv\|_2^2\notag\\ &+C\sum_{p=0}^{P}\sum_{m=2}^{M-1} \|\partial^p_t\nabla^{m-1}v\|^2_2\|\partial^{P-p}_t\nabla^{M+1-m}v\|^2_2 +\frac{\nu}{4}\|\partial^{P}_t\nabla^{M+1}v\|^2_2\notag
\end{align}
In the case $M=0$ the estimate is
\begin{align}
I_{0,P}+J_{0,P}\leq C\sum_{p=0}^{P}\|\partial_t^p v\|_2^2\|\partial_t^{P-p}\nabla v\|_2^2+\frac{\nu}{4}\|\partial^{P}_t\nabla v\|^2_2\notag
\end{align}

We have shown in the case $M>0$
\begin{align}\label{VCHE:inductivedecaystart}
\frac{d}{dt}\|\partial^{P}_t \nabla^{M} v\|^2_2+\frac{\nu}{2}\|\partial^{P}_t
\nabla^{M+1} v\|_2^2 &\leq C\sum_{p=0}^{P}\|\partial_t^p v\|_2^2\|\partial_t^{P-p}\nabla^{M+1} v\|_2^2\\
&+C\sum_{p=0}^P\|\partial_t^p\nabla v\|_2^2\|\partial_t^{P-p}\nabla^Mv\|_2^2 \notag\\
&+C\sum_{p=0}^{P}\sum_{m=2}^{M} \|\partial^p_t\nabla^{m-1}v\|^2_2\|\partial^{P-p}_t\nabla^{M+1-m}v\|^2_2 \notag
\end{align}
and in the case $M=0$,
\begin{align}\label{VCHE:inductivedecaystart0}
\frac{d}{dt}\|\partial^{P}_t v\|^2_2+\frac{\nu}{2}\|\partial^{P}_t\nabla v\|_2^2 \leq C\sum_{p=0}^{P}\|\partial_t^p v\|_2^2\|\partial_t^{P-p}\nabla v\|_2^2
\end{align}

We now begin the inductive part of our argument where the base case is Theorem \ref{VCHE:spacedecay}.  Pick $P\leq K/2$ and assume (inductive assumption) the decay
\begin{equation}\label{inductioninductiondecay}
\|\partial^p_t\nabla^m v\|^2_2\leq C (t+1)^{-2p-m-n/2}
\end{equation}
holds for all $p<P$ and $m$ such that $2p+m\leq K$.  The inductive claim is that the decay holds for $p=P$ with $m$ such that $2P+m\leq K$.  To prove the inductive claim it will be shown first that the decay rate holds for $p=P$ and $m=0$ using (\ref{VCHE:inductivedecaystart0}).  Then, using (\ref{VCHE:inductivedecaystart}) it will be shown that the decay rate holds for the remaining values of $m$ using another inductive argument.

To establish the decay for $p=P$ and $m=0$, apply the inductive assumption to (\ref{VCHE:inductivedecaystart0}) to find
\begin{equation}
\frac{d}{dt}\|\partial^{P}_t v\|^2_2+\frac{\nu}{2}\|\partial^{P}_t\nabla v\|_2^2 \leq C(1+t)^{-n/2}\|\partial_t^P\nabla v\|_2^2 +C(1+t)^{-2P-1-n}\notag
\end{equation}
Take $t$ large enough so that $C(1+t)^{-n/2}\leq \nu/4$ and move the first term on the RHS to the left side
\begin{equation}
\frac{d}{dt}\|\partial^{P}_t v\|^2_2+\frac{\nu}{2}\|\partial^{P}_t\nabla v\|_2^2 \leq C(1+t)^{-2P-1-n}\notag
\end{equation}
Now, an application of Theorem \ref{FSM} with Lemma \ref{VCHE:vtfrequencybound} establishes the decay (\ref{inductioninductiondecay}) for $p=P$ and $m=0$.  This is the base case for the next inductive argument.  Assume (inductive assumption) the decay (\ref{inductioninductiondecay}) holds for $m\leq M+1$ when $p<P$, and $m<M$ when $p=P$, we will show that this implies the decay holds for $m=M$ and $p=P$.  Proving this inductive claim will finish the proof.  Begin by applying the inductive assumption to (\ref{VCHE:inductivedecaystart}).
\begin{align}
\frac{d}{dt}\|\partial^{P}_t \nabla^{M} v\|^2_2+\frac{\nu}{2}\|\partial^{P}_t
\nabla^{M+1} v\|_2^2 &\leq C(1+t)^{-n/2}\|\partial_t^P\nabla^{M+1} v\|_2^2\notag\\
&+C(1+t)^{-n/2-1}\|\partial_t^P\nabla^Mv\|_2^2 \notag\\
&+C(1+t)^{-2p-M-n}\notag
\end{align}
Take $t$ large so that $C(1+t)^{-n/2}\leq \nu/4$ and move the first term on the RHS to the LHS.  Then apply Theorem \ref{FSM} with Lemma \ref{VCHE:vtfrequencybound} to establish the decay rate
\begin{equation}
\|\partial^p_t\nabla^m v\|^2_2\leq C (t+1)^{-n/2}
\end{equation}
Another bootstrap argument gives the optimal decay and finishes the proof.
\end{proof}

\section{Convergence of the VCHE to the NSE in the Whole Space}
To understand how solutions of the VCHE approach solutions of the NSE as the filter constant $\alpha$ approaches zero we must first understand how a solution $u$ of the Helmholtz equation
\begin{equation}\label{Helm:PDE}
u-\alpha^2\triangle u= v
\end{equation}
approaches $v$ when taking $\alpha$ to zero.  To begin we state a theorem concerning the Helmholtz equation in all
of space, the theorem is standard elliptic theory and no proof is
given.  This theorem can be proved using elliptic estimates and
interpolation or if one multiplies the Helmholtz equation by
$e^{\tau/\alpha^2}$ and divides by $\alpha^2$ it can be thought of
as the heat equation and the bounds follow from estimates on the
heat kernel.

\begin{Theorem}\label{Helmholtz:generaltheorem}
Given $v\in L^p(\mathbb{R}^n)$, $p\in (1,\infty)$, there exists a
$u\in W^{1,p}(\mathbb{R}^n)$ that is a weak solution to the
Helmholtz equation $u-\alpha^2\triangle u=v$.  Moreover, this function
satisfies
\begin{align}
\|u\|_p&\leq \|v\|_p \notag\\
\|u\|_q &\leq \frac{C(n,p,q)}{\alpha^{1+\gamma}}\|v\|_p \ for \ \gamma=\frac{n}{2}(\frac{1}{p}-\frac{1}{q})<1 \notag\\
\|\nabla u\|_q &\leq \frac{C(n,p,q)}{\alpha^{3/2+\gamma}} \|v\|_p \
for \ \gamma =\frac{n}{2}(\frac{1}{p}-\frac{1}{q})<\frac{1}{2}\notag
\end{align}
If $n(2/p-1)<1$ then the solution is unique.
\end{Theorem}
\begin{proof}
Standard elliptic theory.
\end{proof}

A solution $u$ of the Helmholtz equation corresponding to $v$ will
approach $v$ weakly as the filter parameter tends to zero. Indeed,
fix $v\in L^p(\mathbb{R}^n)$ and let $\{\alpha_i\}$ be a sequence
tending to zero. By the above theorem, for each $\alpha_i$ there is a
weak solution $u_{\alpha_i}\in W^{1,p}(\mathbb{R}^n)$ of the Helmholtz
equation such that
\begin{equation}
\int_{\mathbb{R}^n} u_{\alpha_i}\cdot \phi \,dx+
\alpha_i^2\int_{\mathbb{R}^n}\nabla u_{\alpha_i}\cdot \nabla \phi \,dx=
\int_{\mathbb{R}^n} v\cdot \phi \,dx\notag
\end{equation}
The functions $u_{\alpha_i}$ are bounded in $L^p(\mathbb{R}^n)$
independent of $\alpha_i$, so there exists a (possible) subsequence $\alpha_{i_j}$
with a weak limit in $L^p(\mathbb{R}^n)$.  Also, for $1/p+1/q=1$
\begin{equation}
\alpha_i^2\int_{\mathbb{R}^n}\nabla u\cdot \nabla \phi \,dx \leq
\alpha_i^2\|\nabla u\|_p\|\nabla \phi\|_q\leq
C(n)\alpha_i^{1/2}\|v\|_p\|\nabla\phi\|_q\notag
\end{equation}
which approaches zero as $\alpha_i\rightarrow 0$.  This proves that
$u_{\alpha_i}\rightharpoonup v$ in $L^p(\mathbb{R}^n)$.  We can
do better then this if $v$ is sufficiently differentiable.

\begin{Theorem}\label{convergencetheorem}
Let $v\in W^{1,p}(\mathbb{R}^n)$ and let $u$ be the corresponding
solution to the Helmholtz equation (\ref{Helm:PDE}).  Then
\begin{equation}
\|u-v\|_q \leq C(n,p,q)\alpha^{1/2-\gamma} \|\nabla v\|_p \ for \
\gamma =\frac{n}{2}(\frac{1}{p}-\frac{1}{q})<\frac{1}{2}\notag
\end{equation}
If $\alpha$ is a sequence tending to zero and $u_\alpha$ are
solutions the Helmholtz equation, then
$u_\alpha\rightarrow v$ strongly in $L^q(\mathbb{R}^n)$ for
$1/p-1/q<1/n$.
\end{Theorem}
\begin{proof}
If $u$ and $v$ satisfy the Helmholtz equation, then
\begin{equation}\label{Helmholtz:quickbound}
\|u-v\|_q\leq \alpha^2\|\triangle u\|_q
\end{equation}
The Helmholtz equation is linear, so the derivatives of the
functions obey the relation $\nabla u -\alpha^2 \triangle\nabla u = \nabla v$.  Applying Theorem \ref{Helmholtz:generaltheorem} to this PDE with the
restriction on $\gamma$ allows the bound
\begin{equation}
\|\triangle u\|_q \leq \frac{C(n,p,q)}{\alpha^{3/2+\gamma}} \|\nabla
v\|_p\notag
\end{equation}
Together with (\ref{Helmholtz:quickbound}),
\begin{equation}
\|u-v\|_q \leq C(n,p,q)\alpha^{1/2-\gamma} \|\nabla v\|_p\notag
\end{equation}
The second statement is an immediate consequence of this.
\end{proof}

In \cite{MR1837927}, \cite{MR1878243}, the authors show how the
solutions of the VCHE approach a solution of the NSE weakly when the
parameter in the filter tends to zero ($\alpha\rightarrow 0$). We
will show how solutions to the VCHE approach solutions to the NSE
strongly as $\alpha\rightarrow 0$ when the solution to the NSE is
known to be regular.  The proof requires estimates on the solution
of the VCHE which are independent of $\alpha$, but in regions of
time where the NSE is known to be regular by some functional
analytic arguments, the passive bound on the filter make this
assumption reasonable.  

For example, solutions of the Navier-Stokes equation obey the Prodi Inequality
 \cite{MR0189354}
\begin{equation}
\frac{d}{dt}\|\nabla u\|_2^2 \leq C\|\nabla u\|_2^{2n} \notag
\end{equation}
This can be used to prove existence of a strong solution in some
time interval $[0,T]$ or regular solutions for all time if the initial data is small.  The Prodi inequality is proved through
energy estimates, using the passive bound for the filter in Theorem
\ref{Helmholtz:generaltheorem} and following the same energy
arguments allows the same bound for solutions of the VCHE.  This
bound will be independent of $\alpha$, so we can apply the following
theorem to conclude that in some closed interval $[0,T]$ the
solution of the VCHE approaches a solution to the NSE strongly.

\begin{Theorem}
Let $\{\alpha_i\}$ be a sequence of filter coefficients tending to zero and $v_{\alpha_i}$ the solutions of the VCHE (\ref{VCHE:PDE}) constructed in Theorem \ref{VCHE:existence} with $\Omega=\mathbb{R}^n$ corresponding to
$w_0\in H^1_\sigma(\mathbb{R}^n)$.
Let $w$ be the solution the NSE with initial conditions $w_0$.
In any time interval $[0,T]$ where a
solution to the NSE is known to be regular, if there exists a bound
\begin{equation}
\sup_{\alpha_i}\sup_{t\in [0,T]}(\|v_{\alpha_i}\|_l+\|\nabla v_{\alpha_i}\|_l) < C \notag
\end{equation}
which is independent of $\alpha$, then $v_\alpha$ approaches $w$ strongly
in $L^\infty([0,T],L^q(\mathbb{R}^n))$ as $\alpha\rightarrow 0$,
where $q=\frac{2l}{l-2}$.
\end{Theorem}
\begin{proof}
We begin with a mild form of the solutions to both problems.  We are working in a time domain with known regularity so these are the
unique solutions.  If $\mathbb{P}$ is the Leray projector onto the divergence free subspace of
$L^2$ and $\Phi$ is the heat kernel, then
\begin{align}
w(t)&=\Phi(t)\ast w_0 -\int_0^t \Phi(t-s)\ast \mathbb{P}\left[w\cdot\nabla w\right](s)\,ds \notag \\
v(t)&=\Phi(t)\ast w_0 -\int_0^t \Phi(t-s)\ast
\mathbb{P}\left[u\cdot\nabla v +\sum u_j\nabla v_j\right](s)\,ds \notag
\end{align}
By adding and subtracting cross terms we see
\begin{align}
w(t)-v(t)= -\int_0^t\Phi(t-s)&\ast \mathbb{P}\left[(w-u)\cdot\nabla w + u\cdot\nabla (w-v)\right] \notag\\
&+ \mathbb{P}\left[u_j\nabla (v_j-w_j) +(u_j-w_j)\nabla w_j \right]
(s) \,ds \notag
\end{align}
The first term in the integrand is bounded using Young's inequality
and the definition of the projector
\begin{align}
\|\Phi(t-s)\ast\mathbb{P}\left[(w-u)\cdot\nabla w\right](s)\|_q
&\leq \|\Phi(t-s)\|_p\|(w-u)\cdot\nabla w\|_2  \notag\\
&\leq\|\Phi(t-s)\|_p\|w-u\|_q\|\nabla w\|_l \notag
\end{align}
where $1/q+1=1/p+1/2$ and $1/2=1/q+1/l$.  Using Theorem
\ref{convergencetheorem} with $\gamma=(1/2-1/q)n/2<1/2$ we obtain
\begin{equation}
\|\Phi(t-s)\ast\mathbb{P}\left[(w-u)\cdot\nabla w\right](s)\|_q \leq
\|\Phi(t-s)\|_p\|\nabla w\|_l\left(\|w-v\|_q +C
\alpha^{1/2-\gamma}\|\nabla v\|_2\right) \notag
\end{equation}
The fourth term can be bounded in essentially the same way. We
approach the second term in a slightly different way, by first
passing the derivative to the heat kernel.  These functions are
smooth functions of the whole space so the projector will commute
with the derivative.
\begin{equation}
\|\Phi(t-s)\ast \mathbb{P}\left[u\cdot\nabla(w-v)\right]\|_q
=\|\nabla\Phi(t-s)\ast \mathbb{P}\left[u\cdot(w-v)\right]\|_q \notag
\end{equation}
Then by using Young's inequality and the definition of the projector
\begin{equation}
\|\Phi(t-s)\ast \mathbb{P}\left[u\cdot\nabla
(w-v)\right]\|_q\leq\|\nabla\Phi(t-s)\|_p \|u\|_l\|w-v\|_q \notag
\end{equation}
To bound the third term, start with the product rule and again pass
the derivative the heat kernel
\begin{align}
\|\Phi(t-s)\ast \mathbb{P}\left[u_j\nabla(v_j-w_j)\right]\|_q&=\|\nabla\Phi(t-s)\ast
\mathbb{P}\left[u_j(w_j-v_j)\right]\|_q\notag\\
&+\|\Phi(t-s)\ast\mathbb{P}\left[\nabla u_j (w_j-v_j)\right]\|_q \notag
\end{align}
Then, using Young's inequality and again the definition of the
projector
\begin{align}
\|\Phi(t-s)\ast &\mathbb{P}\left[u_j\nabla (v_j-w_j)\right]\|_q\notag \\
&\leq \left(\|\nabla\Phi(t-s)\|_p\|u_j\|_l+\|\Phi(t-s)\|_p\|\nabla u_j\|_l\right) \|w_j-v_j\|_q \notag
\end{align}
Putting all of these bounds together and estimating the heat kernel yields
\begin{align}
\|v-w\|_q &\leq A\alpha^{1/2-\gamma}+
B\int_0^t\frac{1}{(t-s)^{\delta}}\|v-w\|_q(s)\,ds \notag\\
A&=C\int_0^t\|\nabla v\|_2\, ds \notag\\
B&=C\sup_{s\in [0,T]}\left(\|\nabla w\|_l+\|u\|_l+\|\nabla
u\|_l\right) \notag
\end{align}
Here, $\delta=1/2+(1-1/p)n/2 <1$ by the assumption $l>n$.
Application of the Gronwall inequality finishes the proof.  For
example, a modified Gronwall inequality \cite{MR2166670} now shows
\begin{align}
\|v-w\|_q&\leq  A\alpha^{1/2-\gamma}\Upsilon(B\Gamma(1-\delta)t^\delta) \notag\\
\Upsilon(z)&=\sum_{n=0}^\infty \frac{z}{\Gamma (n(1-\delta)+1)} \notag
\end{align}
See also \cite{MR1308857}.
Letting $\alpha\rightarrow 0$ we see that $v\rightarrow w$ strongly
in $L^q(\mathbb{R}^3)$.
\end{proof}

\section{Appendix}
Here we construct a weak solution to the VCHE.  Due to the close relation between the VCHE and
the Navier-Stokes equation, our proof is similar to known existence proofs for the NSE.
See, for example, \cite{MR1308857}, \cite{MR673830}, \cite{MR972259}, \cite{MR0050423}, \cite{Leray},\cite{MR1846644}.
First, we construct solutions on any
bounded $\Omega$ with smooth boundary using the Galerkin method, this is where the Stokes
operator is known to be compact thanks to the Poincar\'{e}
Inequality.  Special care is taken to use inequalities which do
not depend on the size of $\Omega$ so we can use these
solutions to prove existence of a weak solution in unbounded
domains. The only step that requires $\Omega$ bounded is in the
compact inclusion used to obtain the strong
convergence necessary to pass limits through the non-linear term.
This problem is overcome by working in the compact support of the
test functions.

To begin we recall a standard but useful elliptic estimate.
\begin{Remark}\label{Helmholtz:goodfilterlemma}
Let $2\leq n\leq 4$ and $\Omega\subset\mathbb{R}^n$ be an open set
with smooth boundary.  If $u\in H^2_\sigma$, and $v\in L^2_\sigma$
satisfy the Helmholtz equation
\begin{equation}
u-\alpha^2 \triangle u=v\notag
\end{equation}
on $\Omega$, then
\begin{align}
\|u\|_n&\leq C \|v\|_2 \label{Helmholtz:goodfilterlemma1}\\
\|\nabla u\|_n&\leq C \|v\|_2 \label{Helmholtz:goodfilterlemma2}\\
\|u\|_2^2+2\alpha^2\|\nabla u\|_2^2+\alpha^4\|\triangle u\|_2^2&=\|v\|_2^2\label{Helmholtz:goodfilterlemma3}
\end{align}
where the constants $C$ depend only on $\alpha$ and $n$.
\end{Remark}

The stationary Stokes equation
\begin{align}\label{Stokes}
\triangle u +\nabla p = v \notag\\
u|_{\partial\Omega} =0\notag
\end{align}
is known to have a solution $u\in H_\sigma^{1}(\Omega)$ for each
$v\in L_\sigma^2(\Omega)$ when $\Omega$ is an open bounded set.  Solving this PDE defines an operator
$L_\sigma^2(\Omega)\rightarrow H_\sigma^1(\Omega)$.  Composing this with the
compact inclusion $H_\sigma^1(\Omega)\rightarrow L_\sigma^2(\Omega)$
gives a compact and self-adjoint operator $L_\sigma^2(\Omega)\rightarrow L_\sigma^2(\Omega)$, which we call the \emph{Stokes operator}.

\begin{Lemma}\label{VCHE:approxexistence}
Let $\Omega\subset \mathbb{R}^n$ be an open bounded set. There exists an orthonormal basis of
$L^2_\sigma(\Omega)$, $\{\omega_j\}_{j=1}^\infty$, where each $\omega_j$ is an eigenfunction of
the Stokes Operator on $\Omega$.  The associated eigenvalues are all
positive real numbers and the eigenvectors are smooth and approach
zero on the boundary. Let $H_m = span\{\omega_1,...,\omega_m\}$ and
let $P_m$ be the orthogonal projection $P_m:L^2_\sigma(\Omega)\rightarrow H_m$.  Given $v_0\in
C^\infty_0(\Omega)$, for each $m$ there is an approximate solution
\begin{equation}
v_m=\sum_{j=1}^m g_{jm}(t)\omega_j \notag
\end{equation}
and
\begin{equation}
u_m= \sum_{j=1}^m\frac{g_{jm}(t)}{1+\alpha^2\lambda_j}\omega_j
\notag
\end{equation}
where $g_{jm}\in C^1([0,T_m])$ for some time $T_m$.
These approximate solutions satisfy the following relations:
\begin{align}\label{VCHE:approxrelation}
<\partial_t v_m, \omega_i>+ <u_m \cdot\nabla v_m, \omega_i> -<\omega_i \cdot\nabla v_m, u_m> &=\nu<\triangle v_m, \omega_i> \\
v_m(0)&=P_m v_0 \notag
\end{align}
\end{Lemma}
\begin{proof}
Owing to spectral theory the Stokes operator (self-adjoint, compact) has a countable number of positive eigenvalues
$\lambda_i$, and associated smooth, divergence-free eigenfunctions
$\omega_i$ which form a basis for $L^2_\sigma(\Omega)$.  These
functions satisfy the relation
\begin{equation}
-\triangle \omega_i=\lambda_i\omega_i \notag
\end{equation}

To determine the scalars $g_{im}$ we construct a system of $m$
ODE's.
\begin{align}
\frac{dg_{im}}{dt}&+ \nu\lambda_i g_{im} \notag \\
&+\sum_{j,k=1}^m\frac{g_{jm}g_{km}}{1+\alpha^2\lambda_k}\left(<\omega_k\cdot\nabla\omega_j,\omega_i>
-<\omega_i\cdot\nabla\omega_j,\omega_k>\right)=0 \notag
\end{align}
Local existence of solutions to ODE's give existence of solutions
$g_{im}$, which are defined for some time interval $[0,T_m]$.
\end{proof}

The bounds in the next lemma will prove that $T_m$ can be bounded
independent of $m$, and in fact $T_m=\infty$ for all $m$.
\begin{Lemma}\label{VCHE:approxbounds}
For $2\leq n\leq 4$, the approximate solutions constructed in Lemma
\ref{VCHE:approxexistence} have the following bounds, which do not
depend on $T$, $\Omega$ or $m$.
\begin{align}
\|v_m\|_{L^\infty([0,T];L_\sigma^2(\Omega))}+\|\nabla v_m\|_{L^2([0,T];L_\sigma^2(\Omega))}
\leq C(n,\alpha,\nu,\|v_0\|_2) \notag\\
\|\partial_t v_m\|_{L^2([0,T];(H^1_\sigma)'(\Omega))} \leq
C(n,\alpha,\nu,\|v_0\|_2) \notag
\end{align}
\end{Lemma}
\begin{proof}
Similar to formal multiplication of the VCHE (\ref{VCHE:PDE}) by
$u$, multiply (\ref{VCHE:approxrelation}) by
$\frac{1}{1+\alpha^2\lambda_i}g_{im}$, sum, then apply Lemma
\ref{KFNSE:bilinear} to see
\begin{align}\label{VCHE:approxenergy}
\|u_m\|_2^2+ \alpha^2\|\nabla u_m\|_2^2 &+2\nu\int_0^T \|\nabla
u_m\|_2^2 \,dt \notag \\
&+2\alpha^2\nu\int_0^T\|\triangle u_m\|_2^2\,dt
=\|u_0\|_2^2+\alpha^2\|\nabla u_0\|_2^2
\end{align}
With the Poincar\'{e} inequality and (\ref{Helmholtz:goodfilterlemma3} this becomes first bound in the theorem.  Using
\ref{Helmholtz:goodfilterlemma1} we deduce
\begin{equation}\label{Helmholtz:integrateun}
\|u_m\|_n^2 +\int_0^\infty \|\nabla u_m\|_n^2 \,dt <
C(n,\alpha,\nu,\|u_0\|_2,\|\nabla u_0\|_2)
\end{equation}

To bound the derivative start with (\ref{VCHE:approxrelation}).  Any $\phi\in H^1_\sigma$ can be written as
a sum of the $\omega_i$ so each approximate solution satisfies
\begin{equation}
<\partial_t v_m, \phi>+ <u_m \cdot\nabla v_m, \phi> -<\phi \cdot\nabla v_m, u_m> =\nu<\triangle v_m, \phi>\notag
\end{equation}
After integration by parts and applying the H\"{o}lder inequality with the Gagliardo-Nirenberg-Sobolev inequality we find
\begin{equation}
|<\partial_t v_m, \phi>| \leq C\|u_m\|_n \|\nabla v_m\|_2 \|\nabla \phi\|_2 +C\|\nabla v_m\|_2 \|\nabla \phi\|_2\notag
\end{equation}
As $\phi$ was chosen arbitrarily we conclude
\begin{equation}
\|\partial_t v_m\|_{(H_\sigma^1)'}\leq C(\|u_m\|_n\|\nabla v_m\|_2+\|\nabla v_m\|_2)\notag
\end{equation}
This, together with (\ref{VCHE:approxenergy}) and (\ref{Helmholtz:integrateun}), proves the second bound in the theorem.
\end{proof}

\begin{Theorem}\label{VCHE:boundedexistence}
Let $\Omega\in\mathbb{R}^n$, $2\leq n\leq 4$ be a bounded set with
smooth boundary and $v_0\in C_0^\infty(\Omega)$.  Then, there exists
a weak solution to the VCHE (\ref{VCHE:PDE}) in the sense of
Definition (\ref{VCHE:weakdefn1}).
\end{Theorem}
\begin{proof}
Thanks to Lemmas \ref{VCHE:approxexistence} and
\ref{VCHE:approxbounds} we only need to prove the convergence of the
approximate solutions. Lemma \ref{VCHE:approxbounds} shows how the
sequence $v_m$ remains bounded, so using a possible subsequence and
the Banach-Alaoglu theorem there exists a function
\begin{align}
v &\in L^\infty ([0,T];L_\sigma^2(\Omega))\cap L^2 ([0,T];H_\sigma^1(\Omega))\notag\\
\partial_t v &\in L^2([0,T];(H_\sigma^1)'(\Omega)) \notag
\end{align}
such that
\begin{align}
v_m\rightharpoonup v \ &in \ L^\infty ([0,T];L^2_\sigma(\Omega))\ weak\ast \label{VCHE:weakconvergence1}\\
v_m\rightharpoonup v \ &in \ L^2 ([0,T];H^1_\sigma(\Omega)) \ weakly \label{VCHE:weakconvergence2}
\end{align}
We will now show that $v$ is a weak solution to the VCHE (\ref{VCHE:PDE}).

By the construction of our approximate solutions and integration by parts, we know for any
basis vector $\omega_j\in L_\sigma^2(\Omega)$ and any smooth scalar
function of time $\phi_j(t)$ such that $\phi_j(T)=0$,
\begin{align}
\int_0^T<v_m,\phi'\omega_j>\,ds &+ \int_0^T<u\cdot\nabla
v,\phi\omega_j>\,ds + \int_0^T<\phi\omega_j\cdot\nabla u,v>\,ds \notag\\
&+\int_0^T<\nabla v,\nabla\phi\omega_j>\,ds =<v_m(0),\phi(0)\omega_j> \notag
\end{align}
The convergence (\ref{VCHE:weakconvergence1}) and (\ref{VCHE:weakconvergence2}) implies
\begin{align}
\int_0^t<v_m,\phi_j'\omega_j>\,ds&\rightarrow\int^t_0<v,\phi_j'\omega_j>\,ds\notag\\
\int_0^t<\nabla
v_m,\phi_j\nabla\omega_j>\,ds&\rightarrow\int_0^t<\nabla
v,\phi_j\nabla\omega_j>\,ds \notag
\end{align}
Also,
\begin{equation}
<v_m(0),\phi_j(0)\omega_j>=<P_m(v_0),\phi_j(0)\omega_j>\rightarrow <v_0,\phi_j(0)\omega_j>
\end{equation}

Passing through the non-linear terms will require strong convergence, so we use the fact that the bounds in
Lemma \ref{VCHE:approxbounds} imply (see \cite{MR972259}, Lemma 8.2) the existence of a possible
subsequence $v_m$ such that
\begin{align}\label{VCHE:strongconvergence}
v_m\rightarrow v \ in \ L^2 ([0,T];L^2(\Omega)) \ strongly
\end{align}
Theorem \ref{Helmholtz:generaltheorem} give the existence of a function $u$
which satisfies
\begin{equation}
u-\alpha^2\triangle u=v \notag
\end{equation}
Similar to \ref{Helmholtz:goodfilterlemma3},
\begin{equation}
\|u_m-u\|_2^2+\alpha^2\|\nabla(u_m-u)\|_2^2+\alpha^4\|\nabla^2(u_m-u)\|_2^2=\|v_m-v\|_2^2\notag
\end{equation}
In particular, applying the Gagliardo-Nirenberg-Sobolev Inequality shows $\|u_m-u\|_n^2\leq C\|v_m-v\|_2^2$.
This, with the strong convergence (\ref{VCHE:strongconvergence}),
shows how $u_m$ approaches $u$ strongly.

We can now prove the convergence of the non-linear terms
\begin{align}
\int_0^T<u_m \cdot\nabla v_m, \phi_j\omega_j>\,ds&\rightarrow\int_0^T<u \cdot\nabla v, \phi_j\omega_j>\,ds \notag\\
\int_0^t<\phi_j\omega_j \cdot\nabla v_m,
u_m>\,ds&\rightarrow\int_0^t<\phi_j\omega_j \cdot\nabla v, u>\,ds
\notag
\end{align}
Indeed, adding and subtracting the cross terms, then using the
H\"{o}lder Inequality, the Gagliardo-Nirenberg-Sobolev Inequality
\begin{equation}
|<u_m \cdot\nabla v_m, \phi\omega_j>-<u \cdot\nabla v, \phi_j\omega_j>| \leq A_1 + B_1\notag
\end{equation}
\begin{align}
A_1&=|<(u_m-u) \cdot\nabla v_m, \phi_j\omega_j>| \notag\\
&\leq\|u_m-u\|_n \|\nabla v_m\|_2\|\phi_j\omega_j\|_\frac{2n}{n-2} \notag\\
&\leq\|v_m-v\|_2 \|\nabla v_m\|_2\|\phi_j\nabla\omega_j\|_2 \notag
\end{align}
Due to the strong convergence (\ref{VCHE:strongconvergence}) the bound
in Lemma \ref{VCHE:approxbounds}, and the H\"{o}lder inequality, we see $\int_0^TA_1\,ds\rightarrow 0$.
Similarly,
\begin{align}
B_1&=|<u \cdot\nabla (v_m-v), \phi_j\omega_j>| \notag\\
&=|<u \cdot\nabla \phi_j\omega_j,(v_m-v)>| \notag\\
&\leq\|u\|_n \|\phi_j\nabla\omega_j\|_\frac{2n}{n-2}\|v_m-v\|_2 \notag\\
&\leq C\|v\|_2 \|\phi_j\nabla \omega_j\|_\frac{2n}{n-2}\|v_m-v\|_2 \notag
\end{align}
Again, owing to
(\ref{VCHE:strongconvergence}), Lemma \ref{VCHE:approxbounds}, and the H\"{o}lder inequality, $\int_0^TB_1\,ds\rightarrow 0$.  Putting this together,
\begin{equation}
\int_0^T<u_m \cdot\nabla v_m,
\phi_j\omega_j>\,ds\rightarrow\int_0^T<u \cdot\nabla v,
\phi_j\omega_j>\,ds \notag
\end{equation}

The remaining non-linear term is handled in a similar way
\begin{equation}
|<\phi_j\omega_j \cdot\nabla v_m, u_m>-<\phi_j\omega_j \cdot\nabla
v, u>| \leq A_2 + B_2 \notag
\end{equation}
\begin{align}
A_2&= |<\phi_j\omega_j \cdot\nabla v_m, (u_m-u)>| \notag\\
&\leq \|\phi_j\omega_j\|_\frac{2n}{n-2} \|\nabla v_m\|_2\|u_m-u\|_n \notag\\
&\leq C\|\phi_j\nabla\omega_j\|_2 \|\nabla v_m\|_2\|v_m-v\|_2 \notag
\end{align}
\begin{align}
B_2&= |<\phi_j\omega_j \cdot\nabla (v_m-v), u>| \notag\\
&= |<\phi_j\omega_j \cdot\nabla u,v_m-v>| \notag\\
&\leq \|\phi_j\omega_j\|_\frac{2n}{n-2} \|v_m-v\|_2\|\nabla u\|_n \notag\\
&\leq C\|\phi_j\nabla\omega_j\|_2 \|v_m-v\|_2\|v\|_2 \notag
\end{align}
Applying
(\ref{VCHE:strongconvergence}) with Lemma \ref{VCHE:approxbounds} and the H\"{o}lder inequality shows
\begin{equation}
\int_0^T<\phi_j\omega_j \cdot\nabla v_m,
u_m>\,ds\rightarrow\int_0^T<\phi_j\omega_j \cdot\nabla v, u>\,ds
\notag
\end{equation}
Since the $\omega_j$ are dense in $L_\sigma^2$ and $\phi_j$ is an
arbitrary smooth function the proof is complete.
\end{proof}

\begin{Corollary}\label{VCHE:boundedexistenceC}
The conclusions of Theorem \ref{VCHE:boundedexistence} hold with the
relaxed hypothesis $v_0\in L^2_\sigma(\Omega)$.
\end{Corollary}
\begin{proof}
Note that all of the bounds attained in Lemma \ref{VCHE:approxbounds} and used in the proof of the previous theorem
depend only on the $L^2$ norm of the initial data.  Let $v^i_{0}\in C^\infty_0(\Omega)$ be a sequence of functions
approaching $v_0$ strongly in $H^1_0$ such that
\begin{equation}\notag
\|v^i_0\|_{H^1_0}\leq \|v_0\|_{H^1_0}
\end{equation}
Such a sequence can be constructed using standard mollifiers and
cutoff functions. Considering each $v^i_0$ as initial data, Theorem
\ref{VCHE:boundedexistence} and its corollary give the existence of
a weak solution $v^i$ in the sense of Definition
\ref{VCHE:weakdefn1}.  Applying (\ref{VCHE:approxbounds}), we see
that these weak solutions satisfy the bounds
\begin{align}
\|v^i\|_{L^\infty([0,T];L_\sigma^2(\Omega))}+\|\nabla v^i\|_{L^2([0,T];L_\sigma^2(\Omega))}
\leq C(n,\alpha,\nu,\|v_0\|_2) \notag\\
\|\partial_t v^i\|_{L^2([0,T];(H^1_\sigma)'(\Omega))} \leq
C(n,\alpha,\nu,\|v_0\|_2) \notag
\end{align}
and for each $\phi\in H^1_\sigma$ the relation
\begin{align}\label{VCHE:H1approx}
\int_0^T<v^i,\partial_t\phi>\,ds &+ \int_0^T<u^i\cdot\nabla
v^i,\phi>\,ds \\
&+ \int_0^T<\phi\cdot\nabla u^i,v^i>\,ds +\int_0^T<\nabla v^i,\nabla \phi>\,ds =<v_0,\phi> \notag
\end{align}
As before, using the Banach-Alaoglu Theorem and extracting a possible
subsequence implies that there exists a function
\begin{align}
v &\in L^\infty ([0,T];L_\sigma^2(\Omega))\cap L^2 ([0,T];H_\sigma^1(\Omega))\notag\\
\partial_t v &\in L^2([0,T];(H_\sigma^1)'(\Omega)) \notag
\end{align}
such that
\begin{align}
v^i\rightharpoonup v \ &in \ L^\infty ([0,T];L^2_\sigma(\Omega))\ weak\ast\notag\\
v^i\rightharpoonup v \ &in \ L^2 ([0,T];H^1_\sigma(\Omega)) \ weakly\notag
\end{align}
Passing the limits through (\ref{VCHE:H1approx}) follows by the same
steps as in the proof of the previous theorem.
\end{proof}

\begin{Theorem}\label{VCHE:unboundedexistence}
Let $v_0\in L^2_\sigma(\mathbb{R}^n)$  Then, there exists a weak
solution in the sense of Definition \ref{VCHE:weakdefn1}, with
initial data $v_0$ in the whole space $\mathbb{R}^n$, $2\leq n \leq
4$.
\end{Theorem}
\begin{proof}
Let $R_i$ be a sequence tending to infinity and $\chi_{R_i}$ a
smooth cutoff function which is equal to $1$ inside the ball of
radius $R_i-\epsilon$ and zero on the boundary of the ball with
radius $R_i$.  The Corollary \ref{VCHE:boundedexistenceC}
now gives existence of a weak solution $v^{R_i}$ on the ball of
radius $R_i$ with initial conditions $v_0\chi_{R_i}$.  Extend
$v^{R_i}$ to all of $\mathbb{R}^n$ by setting it equal to zero
outside the ball of radius $R_i$.  All of the bounds in Lemma
\ref{VCHE:approxbounds} were found independent of the size of
$\Omega$, so here they hold independent of $R_i$. Using the Banach-Alaoglu
Theorem we have the existence of a function
\begin{align}
v &\in L^\infty ([0,T];L_\sigma^2(\mathbb{R}^n))\cap L^2 ([0,T];H_\sigma^1(\mathbb{R}^n))\notag\\
\partial_t v &\in L^2([0,T];(H_\sigma^1)'(\mathbb{R}^n)) \notag
\end{align}
such that
\begin{align}
v^{R_i}\rightharpoonup v \ &in \ L^\infty ([0,T];L^2_\sigma(\mathbb{R}^n))\ weak\ast\\
v^{R_i}\rightharpoonup v \ &in \ L^2 ([0,T];H^1_\sigma(\mathbb{R}^n)) \ weakly
\end{align}
There exists an orthogonal basis $\{\phi_i\}$ for
$L^2([0,T];(\mathbb{R}^n))$ where each function in the basis is
smooth and has compact support in space.  For $R_i$ larger then the
support of $\phi$,  Theorem \ref{VCHE:boundedexistence} with it's corollary show
\begin{align}
\int_0^T<v^{R_i},\partial_t\phi>\,ds &+ \int_0^T<u^{R_i}\cdot\nabla v^{R_i},\phi>\,ds
+ \int_0^T<\phi\cdot\nabla u^{R_i},v^{R_i}>\,ds \notag\\
&+\int_0^T<\nabla v^{R_i},\nabla\phi>\,ds =<v_0,\phi> \notag
\end{align}
The limit $m\rightarrow\infty$ can be passed through the linear
terms just as before. In the (compact) support of each basis function
$\phi_j$, we have the strong convergence to pass the limit through
the non-linear terms.  A diagonal argument shows this
convergence holds as $R_i\rightarrow\infty$.
\end{proof}

In the above existence theorems, the pressure term can be found by
either taking the divergence of the VCHE and solving the Poisson
equation, or using a famous result of de Rham.  See, for example,
\cite{MR1846644}.

\begin{Theorem}\label{VCHE:spaceregularity}
The solutions to the VCHE constructed in Theorems
\ref{VCHE:boundedexistence} and \ref{VCHE:unboundedexistence},
with initial data in $v_0\in H_\sigma^K$, satisfy the bound
\begin{align}
\|\nabla^M v\|_2^2 +\int_0^t\|\nabla^{M+1}v\|_2^2 \leq C(n,\alpha,\nu,\|v\|_{H_0^K})
\end{align}
for all $M\leq K$.
\end{Theorem}
\begin{proof}
We will do the calculations formally and note that these bounds can
be applied to the approximate solutions constructed in Theorem \ref{VCHE:approxexistence}, this proof proceeds by induction.  The inductive assumption is that the following bound holds for all $m<M$.
\begin{equation}
\|\nabla^m v\|^2_2 +\int_0^T \|\nabla^{m+1} v\|_2^2\,dt \leq C\notag
\end{equation}
The base case ($m=0$) is true by Lemma \ref{VCHE:approxbounds}, we
will now show that it holds for $m=M$.
The bound \ref{Helmholtz:goodfilterlemma1} with the
inductive assumption implies
\begin{equation}
\|\nabla^m u\|^2_n +\int_0^T \|\nabla^{m+1} u\|_n^2\,dt \leq C\notag
\end{equation}
Multiply the VCHE (\ref{VCHE:PDE}) by
$\triangle^{M}v$ and integrate by parts to find
\begin{equation}\label{VCHE:regspaceboot}
\frac{1}{2}\frac{d}{dt}\|\nabla^{M} v\|^2_2+\nu\|\nabla^{M+1} v\|_2^2 \leq I_M +J_M
\end{equation}
\begin{align}
I_M&=\sum_{m=0}^{M}{M\choose m}
<\nabla^mu\cdot\nabla\nabla^{M-m}v,\nabla^{M}v>\notag\\
J_M&=\sum_{m=0}^{M}{M\choose m}
<\nabla^Mv\cdot\nabla\nabla^mu,\nabla^{M-m}v>\notag
\end{align}
The two integrals on the RHS are estimated essentially the same way.
The key difference is that in the first one we use the relation
$<u\cdot\nabla v,v>=0$ while in the second we can place an extra
derivative on $u$.

With application of $<u,\nabla v,v>=0$ the first bound becomes
\begin{align}
I_M&=\sum_{m=1}^{M}{M\choose m}<\nabla^mu\cdot\nabla\nabla^{M-m}v,\nabla^{M}v>\notag
\end{align}
H\"{o}lder's inequality, the Sobolev inequality, and Cauchy's
inequality show
\begin{align}
I_M&\leq C\sum_{m=1}^{M} \|\nabla^m u\|_n\|\nabla^{M+1-m}v\|_2\|\nabla^{M}v\|_\frac{2n}{n-2}\notag\\
&\leq C\sum_{m=1}^{M} \|\nabla^m u\|^2_n\|\nabla^{M+1-m}v\|^2_2+\frac{\nu}{4}\|\nabla^{M+1}v\|_2^2\notag
\end{align}
Similarly for the second term
\begin{align}
J_M&\leq C\sum_{m=0}^{M} \|\nabla^M v\|_\frac{2n}{n-2}\|\nabla^{m+1}u\|_n\|\nabla^{M-m}v\|_2\notag\\
&\leq \frac{\nu}{4}\|\nabla^{M+1}v\|_2^2+C\sum_{m=0}^{M} \|\nabla^{m+1} u\|^2_n\|\nabla^{M-m}v\|^2_2\notag
\end{align}
Equation (\ref{VCHE:regspaceboot}) becomes
\begin{equation}
\frac{d}{dt}\|\nabla^{M} v\|^2_2 \leq C\sum_{m=0}^{M} \|\nabla^{m+1} u\|^2_n\|\nabla^{M-m}v\|^2_2\notag
\end{equation}
The Gronwall inequality with application of the inductive assumption finish the proof.
\end{proof}

\begin{Theorem}\label{VCHE:regularity}
The solution to the VCHE constructed in Theorems
\ref{VCHE:boundedexistence} and \ref{VCHE:unboundedexistence},
with initial data in $v_0\in H_\sigma^K$, satisfies the bounds
\begin{align}
\|\partial^p_t\nabla^m v\|_2^2 +\int_0^t\|\partial^p_t \nabla^{m+1}v\|_2^2 \leq C(n,\alpha,\nu,\|v\|_{H_0^K})
\end{align}
for all $M+2P\leq K$.
\end{Theorem}
\begin{proof}
To prove this, we will bound the time derivatives of the solution in terms of the space derivatives, then use the previous theorem to establish regularity.  We will do the calculations formally and note that these bounds can
be applied to the approximate solutions constructed in Theorem \ref{VCHE:approxexistence}.

Apply $\partial_t^P\nabla^M$ to the solution
of the VCHE, from this we have the inequality
\begin{equation}
\|\partial_t^{P+1}\nabla^Mv\|_2^2 \leq
C(\|\partial_t^P\nabla^{M+2}v\|_2^2 +\|\partial_t^P\nabla^M(u\cdot\nabla v)\|_2^2 + \|\partial_t^P\nabla^M(v\cdot\nabla u^T)\|_2^2)\notag
\end{equation}
Using the Gagliardo-Nirenberg-Sobolev inequality and
\ref{Helmholtz:goodfilterlemma} we can bound the first term on the
right hand side as
\begin{align}
\|\partial_t^P\nabla^M(u\cdot\nabla v)\|_2^2
&= \sum_{p=0}^{P}\sum_{m=0}^{M} {P\choose p}{M\choose m}
\|\partial^p_t\nabla^mu\|_n^2 \|\partial^{P-p}_t\nabla^{M+1-m}v\|^2_\frac{2n}{n-2}\notag\\
&\leq C\sum_{p=0}^{P}\sum_{m=0}^{M}
\|\partial^p_t\nabla^mv\|_2^2 \|\partial^{P-p}_t\nabla^{M+2-m}v\|^2_2\notag
\end{align}
Similarly for the second term,
\begin{align}
\|\partial_t^P\nabla^M(v\cdot\nabla u^T)\|_2^2
&= \sum_{p=0}^{P}\sum_{m=0}^{M} {P\choose p}{M\choose m}
\|\partial^p_t\nabla^{m+1}u\|_n^2 \|\partial^{P-p}_t\nabla^{M-m}v\|^2_\frac{2n}{n-2}\notag\\
&\leq C\sum_{p=0}^{P}\sum_{m=0}^{M}
\|\partial^p_t\nabla^{m+1}v\|_2^2 \|\partial^{P-p}_t\nabla^{M+1-m}v\|^2_2\notag
\end{align}
Putting this together we can deduce
\begin{equation}
\|\partial_t^{P+1}\nabla^Mv\|_2^2\leq C\|\partial_t^Pv\|^2_{H^{M+2}_0}\notag
\end{equation}
This implies, for all $M$, $P$, such that $M+2P\leq K$,
\begin{equation}
\|\partial_t^P\nabla^Mv\|_2^2\leq C\|v\|_{H_0^K}^2\notag
\end{equation}
Appealing to Theorem \ref{VCHE:spaceregularity} finishes the proof.
\end{proof}

The previous theorem demonstrates how the norms $\|v\|_{H^m}$ and $\|v\|_{H^{m+1}}$ can be
bounded in terms of $\|v_0\|_{H^m}$.  Since the PDE is parabolic we
can expect regularity from interior estimates but the bounds will
not depend explicitly on the initial conditions.

\begin{Theorem}\label{VCHE:uniequnesstheorem}
The solution to the VCHE constructed in Theorems
\ref{VCHE:boundedexistence} and \ref{VCHE:unboundedexistence} is
unique.
\end{Theorem}

\begin{proof}
Let $v$ and $w$ be two solutions to the VCHE \ref{VCHE:PDE} with the
same initial conditions.  Let $u$ and $\omega$ be the corresponding
``filtered'' velocities.  The difference solves the PDE
\begin{equation}
(v-w)_t-\nu\triangle (v-w) +\nabla p +u\cdot\nabla v
-\omega\cdot\nabla w +v\cdot \nabla u^T -w \cdot \nabla \omega^T =0
\notag
\end{equation}
with zero initial conditions.  Multiplying this relation by $v-w$
and integrating by parts leaves
\begin{align}
\frac{1}{2}\frac{d}{dt}\|v-w\|_2^2+\nu\|\nabla(v-w)\|_2^2 &=
<(u-\omega)\cdot\nabla w, (v-w)>  \notag\\
&+<(v-w)\cdot\nabla (u-\omega),v>\notag\\
&+<(v-w)\cdot\nabla\omega,v-w>
\notag
\end{align}
Using H\"{o}lder's inequality, the
Gagliardo-Nirenberg-Sobolev inequality, Cauchy's inequality, and (\ref{Helmholtz:integrateun}), estimate the RHS
\begin{align}
<(u-\omega)\cdot\nabla w, (v-w)>&\leq \|u-\omega\|_n\|\nabla w\|_2\|v-w\|_\frac{2n}{n-2} \notag\\
&\leq C\|v-w\|_2^2\|\nabla w\|_2^2 +\frac{\nu}{4}\|\nabla(v-w)\|_2^2 \notag\\
<(v-w)\cdot\nabla (u-\omega),v>&\leq \|v-w\|_2\|\nabla(u-\omega)\|_n\|v\|_\frac{2n}{n-2} \notag\\
&\leq C\|v-w\|_2^2\|\nabla v\|_2^2 +\frac{\nu}{8}\|\nabla(v-w)\|_2^2 \notag\\
<(v-w)\cdot\nabla\omega,v-w>&\leq \|v-w\|_2\|\nabla\omega\|_n\|v-w\|_\frac{2n}{n-2} \notag\\
&\leq C\|v-w\|_2^2\|w\|_2^2 +\frac{\nu}{8}\|\nabla(v-w)\|_2^2 \notag
\end{align}
After using the bounds in Lemma \ref{VCHE:approxbounds} we have
\begin{equation}
\frac{1}{2}\frac{d}{dt}\|v-w\|_2^2+\frac{\nu}{2}\|\nabla(v-w)\|_2^2
\leq C\|v-w\|^2_2 \notag
\end{equation}
By assumption $\|v_0-w_0\|_2=0$, so $\|v-w\|_2=0$ for all
$t\in[0,T]$.
\end{proof}

\emph{Acknowledgment}  The authors would like to thank an anonymous referee for helpful advice and comments.

\bibliographystyle{plain}
\bibliography{vche-decay-final}

\end{document}